\documentstyle[11pt]{article}
\newtheorem{thm}{Theorem}[section]
\newtheorem{dfn}{Definition}[section]
\newtheorem{prop}{Proposition}[section]
\newtheorem{lem}{Lemma}[section]
\newtheorem{rem}{Remark}[section]
\newtheorem{prob}{Problem}[section]
\newtheorem{conj}{Conjecture}[section]
\newtheorem{ex}{Example}[section]

\begin{document} 
\title{{\bf Noncommutative algebras related with Schubert 
calculus 
on Coxeter groups}} 
\author{Anatol N. Kirillov and Toshiaki Maeno}
\date{} 
\maketitle 
\begin{abstract}
For any finite Coxeter system $(W,S)$ we construct a certain 
noncommutative algebra, the so-called {\it bracket algebra}, together 
with a family of commuting 
elements, the so-called {\it Dunkl elements.} The Dunkl elements 
conjecturally generate an algebra which is canonically isomorphic to 
the coinvariant algebra of the Coxeter group $W.$ 
We prove this conjecture for classical 
Coxeter groups and $I_2(m)$. We define a ``quantization'' and 
a multiparameter deformation 
of our construction and show 
that for Lie groups of classical type and $G_2,$ the 
algebra generated by Dunkl's elements in the quantized bracket algebra
is canonically isomorphic to 
the small quantum cohomology ring of the corresponding flag variety, 
as described by B. Kim. For crystallographic Coxeter systems 
we define the so-called
{\it quantum Bruhat representation} of the corresponding
bracket algebra. We study in more detail the structure of the relations in
$B_n$-, $D_n$- and $G_2$-bracket algebras,  and as an 
application,
discover {\it a Pieri-type formula} in the $B_n$-bracket algebra. As a
corollary, we obtain  a Pieri-type formula for 
multiplication of an arbitrary $B_n$-Schubert class by some special
ones. Our Pieri-type formula  is a generalization of Pieri's 
formulas obtained by A. Lascoux and M.-P. Sch\"utzenberger for 
flag varieties of type $A.$ 
We also introduce a super-version 
of the bracket algebra together with a family of pairwise 
anticommutative elements, the so-called {\it flat connections with 
constant coefficients}, which describes ``a noncommutative 
differential geometry on a finite Coxeter group'' in a sense 
of S. Majid. 
\end{abstract}
\section*{Introduction}

 \  The study of small quantum cohomology ring of flag varieties 
of type $A$ was initiated by P. Di Francesco and C. Itzykson 
\cite{DFI}, and completed by A. Givental and B. Kim \cite{GK}. 
Later, results of \cite{GK} were generalized by B. Kim \cite{Kim} 
to the case of flag varieties corresponding to any finite 
dimensional semi-simple Lie group. More ``geometric'' approach 
to a description of the small quantum cohomology ring of flag varieties 
was developed in still unpublished lectures by D. Peterson 
\cite{Pe}. Pure algebraic approach to the study of small quantum 
cohomology ring of flag varieties of type $A$ was developed in 
\cite{FGP} and \cite{KM}. A new point of view on both 
classical and quantum cohomology rings of flag varieties of 
type $A$ has been developed in \cite{FK}. 
Namely, the cohomology rings 
in question were realized as certain commutative subalgebras 
in some (noncommutative) quadratic algebras. The latter 
quadratic algebra corresponding to the classical cohomology 
ring of flag variety of type $A,$ 
has many interesting combinatorial and 
algebraic properties, e.g. it appears to be a braided Hopf 
algebra over symmetric group, see e.g. \cite{AS}, \cite{Maj}, \cite{MS}; 
its commutative quotient 
is isomorphic to the algebra of Heaviside's functions 
of hyperplane arrangements of type $A,$ 
see, e.g. \cite{Kir} and the literature quoted therein; the value
of Schubert polynomials on Dunkl elements in the $A_n$-bracket algebra
can be used to describe the structural constants for the product 
of Schubert classes in the cohomology ring of the flag variety 
of type $A_n$ \cite{FK}. 
The main algebraic problem related with the latter 
quadratic algebra is the following: Is this quadratic 
algebra finite dimensional or not?  The main combinatorial problem related
with the bracket algebra $BE(A_n)$ is to find a combinatorial description,
i.e. a ``positive expression'' in the algebra $BE(A_n),$ for the Schubert
polynomials evaluated at the Dunkl elements. It seems natural to raise a 
question: Does there exist for any Coxeter group $W$ a certain algebra 
with properties similar to those for the algebra $BE(A_n)$ ?

In the present paper we are going to present partial answers on the 
questions stated above.
We  introduce and study a generalization 
of the quadratic algebras from \cite{FK} to the case of 
any finite Coxeter system $(W,S).$ Our starting point 
is a remarkable result by C. Dunkl \cite{Du} that the algebra 
generated by ``truncated Dunkl operators'' [{\it ibid\/}] is 
canonically isomorphic to the coinvariant algebra of the 
Coxeter group $W.$ 
It is an attempt to construct a ``quantum coinvariant algebra'' of
a finite Coxeter group and find a ``quantum'' analog of 
C. Dunkl's result mentioned above, that were the main motivation for the
present paper. 

Let us say a few words about the content of our paper.  

In Section 2 we present a definition of bracket algebra $BE(W,S),$ 
as well as that of its super-version $BE^+(W,S),$ 
corresponding to any finite Coxeter system $(W,S).$ If 
$(W,S)$ is the Coxeter system of type $A_n,$ the bracket 
algebra $BE(W,S)$ coincides with the quadratic algebra 
${\cal E}_{n+1}$ introduced and studied in \cite{FK} and 
\cite{Kir}, while the algebra $BE^+(W,S)$ coincides with the quadratic
algebra $\Lambda_{quad}$ of \cite{Maj}, see also \cite{AS}, \cite{MS}. 
We note that the algebra $BE(W,S),$ as well as that $BE^+(W,S),$ is a 
quadratic one only if the Coxeter group $W$ corresponds to a simply-laced 
semi-simple Lie group.
In the case of a crystallographic Coxeter system $(W,S),$ 
except $G_2,$ we introduce
a Hopf algebra structure on the twisted group ring $BE(W,S)\{W \}$ of the
algebra $BE(W,S),$ and show that the latter algebra 
satisfies a ``factorization property'', see Lemma 2.1. As a corollary, 
for a crystallographic Coxeter system $(W,S),$ except $G_2,$ 
we obtain a decomposition of the bracket algebra $BE(W,S)$ into the 
tensor product of certain algebras corresponding to the connected 
components of the Dynkin diagram of the Coxeter system $(W,S)$ 
after removing all the simple edges. 
Our results may be considered as a partial
generalization of results obtained in \cite{FP} 
for $A_n$-quadratic algebras.   

In Section 3 we describe two basic 
representations of the algebra $BE(W,S),$ namely, 
Calogero-Moser's and Bruhat's representations. The latter 
one is a bridge between the algebra $BE(W,S)$ and the Schubert 
Calculus on the Coxeter system $(W,S).$   

 In Section 4 for any $s\in S,$ we introduce Dunkl elements 
in the algebra $BE(W,S),$ denoted by $\theta_s,$ and prove that 
they commute with each other, see Theorem 4.1. The commutative subalgebra 
generated by Dunkl elements is the main object of our study. We also
remark that in the algebra $BE^+(W,S)$ the corresponding elements 
$\theta_s,$ $s\in S,$ are pairwise anticommutative.

In Section 5 we state a ``classical version'' of one of the main 
results of our paper, namely, that for classical Coxeter groups and 
$I_2(m),$ the algebra generated by the
Dunkl elements is canonically isomorphic to the coinvariant 
algebra of the corresponding Coxeter group, see Theorem 5.1. We believe 
that the same result 
holds for any finite Coxeter system. Our proof of Theorem 5.1 is based on
explicit calculations in the corresponding bracket algebras, and we hope 
to improve our techniques
to cover other cases. More specificially, using the defining relations 
in the algebra $BE(B_n),$ we show that all power sums 
$p_{2m}:= \theta_1^{2m}+\cdots +\theta_n^{2m},$ $m>0,$ are equal to zero. 
Note that to show the equality $p_4=0$ in the algebra $BE(B_n),$ 
$n\geq 2,$ we have to use the 4-term 
relations of degree four in the algebra $BE(B_n).$ However, 
in the algebra $BE(D_n),$ $n\geq 4,$ the equality $p_4=0$ follows only from 
quadratic relations. 

In Section 6 we construct a 
quantization $qBE(W,S)$ of our bracket algebra $BE(W,S).$ 

From Section 7 we will assume that Coxeter system $(W,S)$ is a
crystallographic one. Under the assumption that $(W,S)$ is a
crystallographic, we construct 
a representation of the quantized bracket algebra $qBE(W,S)$ 
in the group ring of $W,$ Theorem 7.1. 
The main reason why we made such an assumption on Coxeter system 
$(W,S)$ is that 
the quantum Bruhat representation of the quantized bracket algebra 
$qBE(W,S),$ as defined in Section 7, does not work for general 
noncrystallographic groups, e.g. for $I_2(m),$ if $m\geq 9.$ 
In Section 7 we also state one of the main results of the paper, Theorem 7.2, 
namely, that under the same assumptions as in Theorem 5.1, 
the subalgebra generated by the Dunkl elements 
in the quantized bracket algebra $qBE(W,S)$ is canonically 
isomorphic to the small quantum cohomology ring of the corresponding 
flag variety. 

 In Section 8 we state the ``quantum Chevalley formula'' and 
prove it for classical Lie groups as a corollary of the 
existence of the quantum Bruhat representation and our  
Theorem 7.2. 

 In Section 9 we describe in more detail the bracket 
algebras for Lie groups of type $B_n,$ $D_n$ and $G_2.$ 
In Subsection 9.2 we are going to make use of an algebraic structure of 
relations in the algebra $BE(B_n)$ to the study of the so-called 
{\it Pieri problem} in the Schubert Calculus. Remind that 
Pieri's problem for a finite Coxeter pair $(W,S)$ means to find 
a generalization of the Chevalley formula, see Section 5, for
multiplication of an
arbitrary Schubert class $X_w,$ $w\in W,$ by the Schubert class 
$X_s$ corresponding to a simple reflection $s \in S,$ to the case of 
multiplication of an arbitrary Schubert 
class $X_w$ by the Schubert class $X_u$ corresponding to 
an element $u \in W$ which has a {\it unique} reduced decomposition. 
For the Coxeter group of type $A,$ a solution to 
Pieri's problem is well-known, see e.g. \cite{LS}, \cite{Po}, \cite{Sot}, 
and is given by the so-called {\it Pieri formula.} The latter formula
may be interpreted as an explicit computation 
of the elementary $e_k({\bf X}_m),$ and the complete $h_k({\bf X}_m),$ 
symmetric polynomials in the bracket algebra $BE(A_n)$ after 
the substitution of the variables ${\bf X}_m=(x_1,\ldots, x_m)$
by the $A_n$-Dunkl elements, see e.g. \cite{FK}, \cite{Po}.  
In Subsection 9.2 we give a partial answer on the $B_n$-Pieri 
problem stated above, namely, we give an explicit (if complicated) 
combinatorial formula for the value of the elementary symmetric polynomials 
of an arbitrary degree and the
complete symmetric polynomials of degree at most two 
in the bracket algebra $BE(B_n)$ after the substitution 
of the variables by the $B_n$-Dunkl elements. Let us observe that if 
we specialize all the generators $[i] \in BE(B_n)$ to zero, we 
obtain a $D_n$-analog of Pieri's formula. If we further 
specialize all the generators $\overline{[i,j]} \in BE(B_n)$ 
to zero, we will come to the Pieri rule of type $A_n.$ It is known 
that for Coxeter groups of classical type, the condition 
that an element $u \in W$ has a unique reduced decomposition is 
equivalent to the condition that modulo the ideal generated 
by the fundamental invariant polynomials, the Schubert class $X_u$ 
is equal to either $e_k({\bf X}_m)$ or $h_k({\bf X}_m)$ for some 
$k$ and $m\leq n,$ up to multiplication by some power of 2. 
Let us remark that our Theorem 9.1 describes Pieri's formula in the
algebra $BE(B_n).$ In order to obtain a Pieri-type formula 
in the corresponding (quantum) cohomology ring one has to apply 
the (quantum) Bruhat representation, see Theorems 3.2 and 7.1. 
Since both the classical and the quantum Bruhat representations have 
a huge kernel, it is not obvious how to deduce the Pieri-type 
formulas of \cite{BS} and \cite{PR} from the $B_n$-type Pieri 
formulas of this paper.  \\
It seems very interesting problems to extend our results 
to the cases of the Grothendieck ring and (quantum) equivariant 
cohomology ring of flag varieties. 
We will consider these problems in subsequent publications. 

We expect that for simply-laced Coxeter systems $(W,S)$ the algebra 
$BE(W,S)$ is a finite dimensional braided Hopf algebra over 
$W.$ However, our algebra $BE(D_4)$ is different from the pointed 
Hopf algebra over $D_4$ constructed in \cite{MS}. Surprisingly, 
the latter Hopf algebra appears to be isomorphic to 
a certain quotient of the algebra $BE^+(B_2),$ see Section 9.1. 
For nonsimply-laced Coxeter systems $(W,S)$ the algebra 
$BE(W,S)$ turns out to be infinite dimensional, but it seems 
plausible that a certain finite dimensional quotient of the algebra 
$BE(W,S)$ has a natural structure of a pointed Hopf algebra, 
and the algebra generated by the images of Dunkl's elements 
is isomorphic to that in the algebra $BE(W,S).$ 

The main motivation for introducing our bracket algebra and its 
quantization is an intimate connection of the former and latter 
with classical and 
quantum Schubert Calculi for finite Coxeter groups \cite{BGG}, 
\cite{Hi}. For Coxeter systems of type $A,$ combinatorial and 
algebraic study of Schubert polynomials was initiated and 
developed 
in great details by Alain Lascoux and Marcel-Paul Sch\"utzenberger 
\cite{LS}. 
It is our pleasure to express deep gratitude to Alain Lascoux from 
whom we have learned a lot about this beautiful and deep branch 
of Mathematics. 
\section{Coxeter groups} 
Most part of this section can be found in Humphreys \cite{Hu}. 
\begin{dfn} A Coxeter system is a pair 
$(W,S)$ of a group $W$ and a set of generators 
$S\subset W,$ subject to relations 
\[ (ss')^{m(s,s')}=1, \] 
where $m(s,s)=1$ and $m(s,s')=m(s',s)\geq 2$ for 
$s\not= s' \in S.$ The group $W$ is called a Coxeter 
group. 
\end{dfn} 
\begin{dfn} Let $(W,S)$ be a Coxeter system. 
For an element $w\in W,$ the number 
\[ l(w)=\min\{r \mid w=s_1\cdots s_r, \; \; s_i \in S\} \] 
is called the length of w. We say the expression 
$w=s_1\cdots s_r$ $(s_i\in S)$ is reduced if $r=l(w).$ 
The set of all reduced expressions of an element $w \in W$ 
is denoted by $R(w).$ 
\end{dfn} 
We assume $S$ to be finite. 
Let $V$ be an ${\bf R}$-vector space with a basis 
$\Sigma=\{ \alpha_s \mid s\in S\}$ and symmetric bilinear 
form $( \; , \; )$ such that 
\[ (\alpha_s,\alpha_{s'})=-\cos \frac{\pi}{m(s,s')}. \] 
Consider the linear action $\sigma$ of 
$W$ on $V$ defined by 
\[ \sigma (s) \lambda = 
\lambda - 2(\alpha_s,\lambda)\alpha_s . \] 
The representation $\sigma : W\rightarrow {\rm GL}(V)$ 
is called the geometric representation of $W.$
\begin{dfn}
We define the root system $\Delta$ 
of $W$ to be the set of the all images of $\alpha_s$ under the 
action of $W.$ 
\end{dfn} 
Any element $\gamma\in \Delta$ can be expressed in the form 
\[ \gamma =\sum_{s\in S}c_s \alpha_s \; (c_s\in {\bf R}). \]
Call $\gamma$ positive (resp. negative) and write 
$\gamma > 0$ (resp. $\gamma < 0$) if all $c_s\geq 0$ (resp. 
$c_s\leq 0$). Write $\Delta_{+}$ (resp. $\Delta_{-}$) for 
the set of positive (resp. negative) roots. Note that 
$\Delta=-\Delta$ and $\Delta =\Delta_{+}\amalg\Delta_{-}.$ 
\begin{lem} 
The representation 
$\sigma : W\rightarrow {\rm GL}(V)$ is faithful. 
\end{lem} 
For a given root $\gamma = w(\alpha_s)$ ($w\in W,$ $s\in S$), 
the element $wsw^{-1}$ depends only on $\gamma$ and it acts 
on $V$ as a reflection sending $\gamma$ to $-\gamma.$ 
We denote it by $s_{\gamma}.$ 
\begin{lem}
Let $w\in W$ and $\gamma \in \Delta_{+}.$ 
Then $l(ws_{\gamma})>l(w)$ if and only if $w(\gamma)>0.$  
\end{lem}
\begin{dfn} 
The Coxeter system $(W,S)$ is called crystallographic when 
its root system $\Delta$ can be normalized to satisfy the 
condition 
\[ \frac{2(\gamma,\gamma')}{(\gamma,\gamma)}\in {\bf Z} \] 
for all $\gamma, \gamma' \in \Delta. $
\end{dfn} 
In our paper, the crystallographic systems are always normalized to 
satisfy the condition above. 
\section{Bracket algebra of Coxeter group} 
\subsection{Definition of the bracket algebra}
\begin{dfn} 
Let $(W,S)$ be a Coxeter system and assume $W$ to be finite. 
We define the bracket algebra $BE(W,S)$ as an 
associative algebra over ${\bf R}$ with generators $[\gamma],$ 
$\gamma\in \Delta,$ subject to the following relations: 
\medskip \\ 
{\rm (i)} For any $\gamma \in \Delta,$ 
\[ [-\gamma]=-[\gamma]. \] 
{\rm (ii)} For any $\gamma \in \Delta,$ 
\begin{equation}
[\gamma]^2=0.
\label{eq0}
\end{equation}
{\rm (iii) (Quadratic relations)} 
Let $\Delta'=\{ \gamma_0, \ldots, \gamma_{m-1} \} 
\subset \Delta_{+}$ be a set of positive roots 
such that ${\bf R}_{\geq 0}\langle \gamma_i,\gamma_{i+1}
\rangle \cap \Delta_{+}= \{ \gamma_i,\gamma_{i+1} \} $ for 
all $i=0,\ldots ,m-2.$ 
If $\Delta'$ forms a root system of type $I_2(m)$ $(m\geq 2),$ 
then 
\begin{equation}
\sum_{i=0}^{m} [\gamma_i][\gamma_{i+k}]=0
\label{eq1} 
\end{equation}  
for $1\leq k \leq m/2,$ where we set by definition 
$\gamma_{j+m}=-\gamma_j.$ \\
{\rm (iv) (4-term relations for subsystems of type $I_2(m)$)} 
Let $\Delta'\subset \Delta_{+}$ be a set of 
positive roots as in {\rm (iii).} 
If $\Delta'$ forms a root system of type $I_2(m),$ $m\geq 4,$ 
and $k=[m/2]-1,$  
then 
\[
[\gamma_k] \cdot [\gamma_0][\gamma_1]\cdots [\gamma_{2k}]
+ [\gamma_0][\gamma_1] \cdots 
[\gamma_{2k}] \cdot [\gamma_k] \]
\[ + [\gamma_k]\cdot [\gamma_{2k}][\gamma_{2k-1}] 
\cdots [\gamma_0]
+ [\gamma_{2k}][\gamma_{2k-1}] \cdots [\gamma_0]\cdot [\gamma_k]
=0. \] 
\end{dfn}
\begin{rem}
{\rm 1) The defining ideal generated by the 
relations (i), (ii), (iii) and (iv) 
is stable with respect to the action of the Weyl group $W.$ 
In other words, the algebra $BE(W,S)$ is a $W$-module. 
\smallskip \\ 
2) If $(W,S)$ is a Coxeter system of type $A_n,$ then the bracket 
algebra $BE(W,S)$ coincides with the quadratic algebra 
${\cal E}_{n+1}$ 
introduced in \cite{FK}, see also \cite{Kir}. \smallskip \\ 
3) All the defining relations above come from the subsystems 
of rank two. As for explicit descriptions of these relations 
in the case of type $B_2$, $D_2$ and $G_2,$ as well as for $B_n$
and $D_n$ types, see Section 9. 
\smallskip \\
4) Algebra $BE(W,S)$ has a natural grading, 
if we consider the generators $[\gamma]$ as elements of 
degree one.} 
\end{rem} 
\begin{prob}{\rm 
Find the Hilbert series of the bracket algebra 
$BE(W,S).$} 
\end{prob} 
We expect that the algebra $BE(W,S)$ is finite dimensional for 
simply-laced Coxeter groups. 
\begin{prob} 
{\rm Describe the algebra $BE(W,S)$ as a $W$-module, find its character, 
and / or the graded multiplicities of 
its irreducible components.} 
\end{prob}
\begin{rem}
{\rm We can define the super-version $BE^+(W,S)$ 
of the bracket algebra by using the relation 
$[\gamma]=[-\gamma]$ 
($\gamma \in \Delta$) instead of (i) in Definition 2.1. If 
$(W,S)$ is of type $A_n,$ the algebra $BE^+(W,S)$ coincides 
with the algebra $\Lambda_{quad}$ of \cite{Maj}, see also 
\cite{AS} and \cite{MS}. For crystallographic groups, one can show 
that the left-invariant Woronowicz exterior algebra $\Lambda_w$ 
\cite{Wo} for some special choice of a differential structure on 
$W,$ see \cite{Maj}, is a quotient of the algebra $BE^+(W,S).$ 
However, in a nonsimply-laced case, the algebra $\Lambda_w$ is a proper 
quotient of our algebra $BE^+(W,S).$ } 
\end{rem} 
\subsection{Hopf algebra structure on the twisted group algebra}
Since the bracket algebra $BE(W,S)$ has a 
$W$-module structure, one can construct the twisted group algebra 
$BE(W,S)\{ W \} = \{ \sum_{w\in W} c_w \cdot w \; 
| \; c_w \in BE(W,S) \} $ 
by putting commutation relations 
$w[\gamma]=[ w( \gamma ) ]w$ for $w\in W$ and 
$[ \gamma ] \in BE(W,S).$ 
\begin{prop}
Let $(W,S)$ be a crystallographic Coxeter system, except $G_2,$ 
the twisted group algebra $BE(W,S)\{ W \}$ 
has a natural Hopf algebra structure with the coproduct $\Delta,$ 
the antipode $S$ and the counit $\epsilon$ 
defined by the following formulas: \smallskip \\ 
\[ 
\begin{array}{ll}
\Delta([ \gamma ]) = [ \gamma ] \otimes 1 + s_{\gamma} 
\otimes [ \gamma ], & \Delta (w) = w \otimes w , \\
S( [ \gamma ] ) = s_{\gamma} [ \gamma ], & S(w)= w^{-1}, \\  
\epsilon( [ \gamma ] ) =0, & \epsilon (w)=1, 
\end{array}
\]
for $[ \gamma ] \in BE(W,S)$ and $w\in W.$
\end{prop} 
Such a Hopf algebra structure was invented and studied in 
\cite{FP} 
for $A_n$-quadratic algebras. 

The Hopf algebra $BE(W,S)\{ W \}$ acts on itself by the adjoint action 
\[ w:x \mapsto wxw^{-1} , \; \; w \in W, \] 
\[ [\gamma]: x \mapsto [\gamma]x-s_{\gamma}(x)[\gamma]. \] 
The subalgebra $BE(W,S)$ is invariant under the adjoint action 
of $BE(W,S)\{ W \}.$ 
The element $[\gamma]\in BE(W,S)$ acts on $BE(W,S)$ by a twisted 
derivation 
\[ D_{\gamma} (x) = [\gamma]x-s_{\gamma}(x)[\gamma], \] 
which satisfies the twisted Leibniz rule 
\[ D_{\gamma} (xy)= D_{\gamma}(x)y+s_{\gamma}(x)D_{\gamma}(y). \]
\begin{lem}
Let $(W',S')$ be a parabolic subsystem of $(W,S)$ and $\Delta'$ the set of
roots corresponding to $(W',S').$ 
Denote by ${\cal A}(\Delta \setminus \Delta ')$ 
the subalgebra of $BE(W,S)$ 
generated by the elements $[\gamma],$ 
$\gamma \in \Delta \setminus \Delta'.$ 
Assume that $S\setminus S' = \{ t \}$ and $m(s,t)\leq 3$ 
for any $s\in S'.$ Then the subalgebra 
${\cal A}(\Delta \setminus \Delta ')$ is invariant under
the adjoint action of algebra $BE(W',S')$, 
and the multiplication map  
\[ [\gamma''] \otimes [\gamma'] \mapsto [\gamma''][\gamma'], \; \; \; 
[\gamma'] \in  BE(W',S'), \; \;
[\gamma''] \in  {\cal A}(\Delta \setminus \Delta ')  \]
defines a $BE(W',S')$-linear isomorphism of 
algebras 
\[ {\cal A}(\Delta \setminus \Delta ') 
\otimes BE(W',S') \cong BE(W,S),  \] 
where $BE(W',S')$-module structure on the tensor product 
${\cal A}(\Delta \setminus \Delta ') \otimes BE(W',S')$ 
is given by 
\[ [ \gamma] (a \otimes b) = D_{\gamma}(a)\otimes b + 
s_{\gamma} \otimes [\gamma]b . \]
\end{lem}
It follows from Lemma 2.1 that the Hilbert series of algebra
$BE(W',S')$ divides that of algebra $BE(W,S).$ We give a few more examples
of application of Lemma 2.1 in Section 9.
\begin{rem}
{\rm It is not difficult to see that the algebra $BE(W,S)$ is a braided 
group in the category of $W$-crossed modules with braiding 
$\Psi ([\gamma_1]\otimes [\gamma_2]) = 
s_{\gamma_1}[\gamma_2] \otimes [\gamma_1].$ The Hopf algebra 
$BE(W,S)\{ W \}$ is obtained as its biproduct bosonization 
in the sense of Majid. For Coxeter groups of type $A$ these results
have been shown originally by S.Majid, see \cite{Maj} and the literature
quoted therein. }
\end{rem}
\section{Representations of bracket algebra} 
In this section we are going to construct two basic 
representations 
of the algebra $BE(W,S).$ 
\subsection{Calogero-Moser representation}
Given the geometric representation 
$\sigma : W\rightarrow {\rm GL}(V),$ it induces the natural 
action of $W$ on the ring of polynomial 
functions ${\bf S}(V^{\ast}).$ For any positive root $\gamma,$ 
the divided difference operator $\partial_{\gamma}$, or 
Demazure's operator \cite {De}, acting 
on the ring ${\bf S}(V^{\ast})$ is defined by 
\[ \partial_{\gamma} = \frac{1-s_{\gamma}}{\gamma}. \] 
\begin{thm}
A map $[\gamma] \mapsto \partial_{\gamma}$ defines 
a representation of the algebra $BE(W,S).$ 
\end{thm} 
{\it Proof.} Compatibility with the relation (iv) is clear. 
As for the compatibility with the relation (iii), we may restrict 
our consideration to subsystems of rank two. It is easy to check 
the compatibility for $A_2,$ $B_2$ and $I_2(m).$ \rule{3mm}{3mm} 
\subsection{Bruhat representation}
Let us define a linear operator ${\bf s}_{\gamma}$ acting 
on the group ring ${\bf R}\langle W \rangle $ by the rule 
\[ {\bf s}_{\gamma}.w = \left\{ 
\begin{array}{ll} 
ws_{\gamma}, & {\rm if} \; \; l(ws_{\gamma})=l(w)+1, \\ 
0, & {\rm otherwise.} 
\end{array} 
\right. \] 
\begin{thm} 
A map $[\gamma] \mapsto {\bf s}_{\gamma}$ defines 
a representation of the algebra $BE(W,S).$  
\end{thm} 
{\it Proof.} It is enough to show the compatibility with 
the relation (iii). We use only linear relations among the 
roots in the subsystem of rank two containing 
$\alpha$ and $\beta.$ 
We may assume that 
$\alpha$ and $\beta$ generate a root system of type $I_2(m)$ 
($m\geq 2$). 
Let $a_i=(\cos(i\pi /m),\sin(i\pi/m))\in {\bf R}^2,$ 
$i=0,\ldots ,m-1.$ Then 
$\Delta_{+}= \{ a_0,\ldots ,a_{m-1} \} $ and 
we have to check 
\[ \sum_{i=m-k}^{m-1}[a_{i+k-m}][a_i]w 
=\sum_{i=0}^{m-k-1} [a_{i+k}][a_i]w , \] 
for $1\leq k \leq (m-1)/2.$ From now on, we put $k=1$ 
for simplicity, but the following argument works well 
for all $k.$ 
If ${\bf s}_{a_{i+1}}{\bf s}_{a_i}w=ws_{a_i}
s_{a_{i+1}},$ then $l(w)=l(ws_{a_i})-1$ and 
$l(ws_{a_i})=l(ws_{a_i}s_{a_{i+1}})-1.$ From Lemma 1.2, 
$w(a_i)>0$ and 
$ws_{a_i}(a_{i+1})=-w(a_{i-1})>0.$ 
So we have that $w(a_{m-1})>0$ and 
$ws_{a_{m-1}}(a_0)=w(a_{m-2})<0,$ and that 
$w(a_j)$ and $w(a_{j-1})$ are both positive or both 
negative for $j\not= i.$ 
Hence, if 
${\bf s}_{a_{i+1}}{\bf s}_{a_i}w=
ws_{a_i}s_{a_{i+1}},$ 
then 
${\bf s}_{a_{j+1}}{\bf s}_{a_j}w=0$ for 
$j\not= i$ and 
${\bf s}_{a_0}{\bf s}_{a_{m-1}}w= 
ws_{a_{m-1}}s_{a_0}.$ Conversely, if 
${\bf s}_{a_0}{\bf s}_{a_{m-1}}w= 
ws_{a_{m-1}}s_{a_0},$ then there is only one $i$ such 
that 
${\bf s}_{a_{i+1}}{\bf s}_{a_i}w=
ws_{a_i}s_{a_{i+1}}$ and 
${\bf s}_{a_{j+1}}{\bf s}_{a_j}w=0$ for $j\not=i.$ 
\rule{3mm}{3mm} 
\begin{prob} 
{\rm Does there exist a finite dimensional faithful 
representation of the algebra  $BE(W,S) ? $} 
\end{prob} 
\section{Chevalley and Dunkl elements} 
\begin{dfn} 
For each $s\in S,$ the Chevalley 
element $\eta_s$ in the algebra $BE(W,S)$ is defined by 
\begin{equation} 
\eta_s = \sum_{\gamma \in \Delta_{+}} 
\langle \omega_s, \gamma^{\vee} \rangle [\gamma] 
\label{eq3} , 
\end{equation} 
where $\omega_s$ is the fundamental dominant weight 
corresponding to $\alpha_s$ and 
$\gamma^{\vee}=2\gamma/(\gamma,\gamma).$ 
\end{dfn}
\begin{dfn}
For each $s\in S,$ the Dunkl 
element $\theta_s$ in the algebra $BE(W,S)$ is defined by 
\[ \theta_s = \sum_{s'\in S}c_{s,s'}\eta_{s'} , \]
where the coefficients $c_{s,s'}$ are defined by 
$c_{s,s'}=(\alpha_s, \alpha_{s'}).$  
\end{dfn} 
\begin{thm}
The Dunkl elements $\theta_s$ $(s\in S)$ 
commute pairwise. 
\end{thm} 
{\it Proof.} It is enough to show that the Chevalley elements 
commute pairwise. 
First of all, let us observe that the element 
$\eta_s\eta_s'-\eta_s'\eta_s$ can be decomposed as a sum of 
contributions from root subsystems of rank two. Thus, 
we may assume that the root system $\Delta$ 
is of type $I_2(m).$ Let $S=\{ a_0 , a_{m-1} \}$ and 
\[ a_i = \lambda_1^{-1}\lambda_{i+1}a_0+ \lambda_1^{-1}
\lambda_i a_{m-1}, \; \; \; \lambda_i=\sin \frac{i}{m}\pi,
\; 0\leq i \leq m-1. \] 
Then $\Delta = \{a_0,a_1,\ldots,a_{m-1} \} .$ We have to 
show that $\eta_1$ and $\eta_2$ commute, where 
\[ \eta_1= \sum_{i=0}^{m-1} \lambda_{i+1}[a_i], \; \; 
\eta_2 = \sum_{i=0}^{m-1} \lambda_i [a_i]. \]
We have 
\[ 2(\eta_1 \eta_2 -\eta_2 \eta_1) = 
\sum_{i,j=0}^{m-1} 
\left( \cos \frac{i-j+1}{m}\pi + \cos \frac{i+j+1}{m}\pi \right)
\left( [a_i][a_j]-[a_j][a_i] \right) .\] 
Here, $\cos((i+j+1)\pi/m)$ is symmetric on $i$ and $j,$ so 
\[ \sum_{i,j} \left( \cos \frac{i+j+1}{m}\pi \right) 
\left( [a_i][a_j]-[a_j][a_i] \right) =0. \] 
Note that $s_{a_i}s_{a_j}=s_{a_p}s_{a_q}$ if and only if 
$i-j\equiv p-q \; {\rm mod} \; m.$ Hence the relations in 
Definition 2.1 (iii) imply that 
\[ \sum_k \sum_{i-j\equiv k \; (m)}
\left( \cos \frac{k+1}{m}\pi \right)
\left( [a_i][a_j]-[a_j][a_i] \right) =0. 
\; \; \; \; \rule{3mm}{3mm} \] 
\begin{rem} 
{\rm For the commutativity of the Dunkl elements 
$\theta_s$ it is enough to assume the validity of quadratic relations 
(iii) in Definition 2.1 only.} 
\end{rem} 
\begin{rem} 
{\rm In a similar fashion one can check that the elements $\theta_s,$ 
$s\in S,$ defined as in Definition 4.2 in the super-version 
$BE^+(W,S)$ of the bracket algebra $BE(W,S)$ are pairwise 
{\it anticommutative.} It is a challenging problem to describe the 
subalgebra in $BE^+(W,S)$ generated by the elements $\theta_s, s\in S.$}
\end{rem}  
\section{Algebra generated by Dunkl elements}
Let $ | S | =n.$ 
In case when $W$ is a finite reflection group, it is known 
that the subalgebra 
${\bf S}(V^{\ast})^{W}\subset{\bf S}(V^{\ast})$
of $W$-invariant polynomials is generated over ${\bf R}$ 
by $n$ homogeneous, 
algebraically independent polynomials $f_1,\ldots ,f_n$ of 
positive degree. We denote by $I_W\subset {\bf S}(V^{\ast})$ 
the ideal generated by $f_1,\ldots,f_n.$ The quotient ring 
${\bf S}_W:={\bf S}(V^{\ast})/I_W$ is called the coinvariant 
algebra of $W.$ 

An explicit construction of a linear basis of ${\bf S}_W$ is 
given by Bernstein, Gelfand and Gelfand \cite{BGG}, and Hiller 
\cite{Hi}. Let $w=s_{i_1}\cdots s_{i_l}$ $(s_{i_1},\ldots ,
s_{i_l}\in S)$ be a reduced decomposition of $w\in W.$ 
We define the operator $\partial_w$ acting on the algebra of polynomial
functions ${\bf S}(V^{\ast})$ 
by $\partial_w =\partial_{\alpha_{s_{i_1}}} \! \cdots
\partial_{\alpha_{s_{i_l}}},$ where 
$\partial_{\alpha_{s_{i_1}}},\ldots ,
\partial_{\alpha_{s_{i_l}}}$ are divided difference operators 
defined in Section 3. The definition of the operator 
$\partial_w$ is independent of the choice of 
a reduced decomposition of $w.$ 

For any polynomial $f \in {\bf S}(V^{\ast}),$ one can define an element
$[f] \in BE(W,S)$ as an image of $f$ by the algebra homomorphism 
obtained by the substitution $\omega_s \mapsto \eta_s.$ 
\begin{dfn} 
{\rm (cf. \cite{BGG},\cite{Hi})} 
We define the polynomials $X_w \in {\bf S}(V^{\ast}),$ 
$w\in W,$ by the following 
formulas: 
\[ X_{w_0}=|W|^{-1}\prod_{\gamma\in \Delta_{+}}\gamma, \] 
\[ X_w=\partial_{w^{-1}w_0}X_{w_0}, \] 
where $w_0\in W $ is the element of maximal length. 
\end{dfn} 
It is known (\cite{BGG}, \cite{Hi}) that 
the images of the polynomials $\{ X_w \}$ 
in the coinvariant algebra ${\bf S}_W$ form a linear basis and satisfy
the {\it Chevalley formula} 
\[ X_s X_w = \sum (\omega_s,\gamma^{\vee}) X_{ws_{\gamma}} 
\; \; \; {\rm mod} \; I_W, \] 
where the sum is taken over the positive roots $\gamma$ such that 
$l(ws_{\gamma})=l(w)+1.$ It is useful to note that one can 
obtain the Chevalley formula above by applying the Bruhat representation, 
see Theorem 3.2, to the equality $[X_s]= \eta_s$ in the algebra 
$BE(W,S).$ 

We have the following statement from the Chevalley formula.  
\begin{lem}
There exists a surjective homomorphism from ${\bf S}_W$ 
to the subalgebra generated by the Chevalley elements 
${\bf R}[\eta_s | s\in S ] \subset BE(W,S),$ 
which maps $X_s$ to $\eta_s.$ 
\end{lem}

\begin{thm} 
For Coxeter groups of classical type and $I_2(m)$, the subalgebra
${\bf R}[\theta_s | s\in S ]$ in $BE(W,S)$ generated 
by Dunkl elements is canonically isomorphic to 
the coinvariant algebra of the group $W,$ 
i.e. 
\[ {\bf R}[ \theta_s \mid s\in S ] 
\tilde{=} {\bf S}_W. \] 
\end{thm} 
We postpone a proof till Section 9.
\begin{conj}
The statement of Theorem 5.1 is valid for any finite Coxeter group.
\end{conj}
\begin{conj}
{\rm Let $(W,S)$ be a crystallographic Coxeter system, then
there exists a monomial basis 
$\{ b_{\mu} \} _{\mu}$ in the algebra $BE(W,S),$ 
such that for any $w\in W$ the polynomial 
$[X_w]$ can be expressed as a linear 
combination of $b_{\mu}$'s with {\it nonnegative} coefficients.}
\end{conj}
\section{Quantization of bracket algebra} 
We consider the group of characters 
\[ C= {\rm Hom}(V^{\ast}, S^1), \] 
and its elements $q_s=\exp 
(2\pi \sqrt{-1}\langle \; \cdot \; ,\alpha_s^{\vee} \rangle)$ 
for $s \in S.$ For $\gamma^{\vee}=
\sum_{s\in S} n_s \alpha_s^{\vee},$ we set $q_{\gamma^{\vee}}=
\prod_s q_s^{n_s}.$ 
\begin{dfn} 
The quantized bracket algebra $qBE(W,S)$ 
is the associative algebra over the ring
${\bf R}[q_s \mid s\in S]$ with generators 
$[\gamma],$ $\gamma \in \Delta,$ subject to the relations: 
\\  
${\rm (i)'}$ For any $\gamma \in \Delta,$ 
\[ [-\gamma]=-[\gamma]. \] 
${\rm (ii)'}$ For $\gamma \in \Delta_+,$ 
\[ 
[\gamma]^2 = q_{\gamma}, \; \; \; {\rm if} \; \; \; 
\gamma \in \Sigma, \] 
\[ [\gamma]^2 =0, \; \; \; {\rm otherwise.} \] 
${\rm (iii)'}$ The same relations as in Definition 2.1 (iii). \\ 
${\rm (iv)'}$ Under the same assumptions as in Definition 2.1 (iv),
if in addition the following inequality \\
$l(s_{\gamma_k})\not=2(\rho,\gamma_k^{\vee})-1 \; \; \; holds, \; \; then$ 
\[  
[\gamma_k] \cdot [\gamma_0] [\gamma_1]\cdots 
[\gamma_{2k}]+ [\gamma_0] [\gamma_1]\cdots 
[\gamma_{2k}] \cdot [\gamma_k] \] 
\[ + [\gamma_k] \cdot [\gamma_{2k}][\gamma_{2k-1}] 
\cdots [\gamma_0]
+[\gamma_{2k}][\gamma_{2k-1}] \cdots [\gamma_0]
\cdot [\gamma_k]
=0. \]     
\end{dfn} 
  
\begin{dfn} 
The Chevalley elements $\tilde{\eta}_s$ and the Dunkl 
elements $\tilde{\theta}_s,$ 
$s\in S,$ in the algebra $qBE(W,S)$ are defined 
by the same formulas as in Definitions 4.1 and 4.2. 
\end{dfn} 
\begin{thm} 
The Dunkl elements $\tilde{\theta}_s,$ 
$s\in S,$ commute pairwise. 
\end{thm} 
{\it Proof.} The proof can be done in the same manner as that of
Theorem 3.1. 
\rule{3mm}{3mm}
\begin{rem} 
{\rm It is natural to consider a multiparameter 
deformation of the algebra $BE(W,S)$ which is generated by 
elements $[\gamma],$ $\gamma \in \Delta,$ with defining 
relations ${\rm (i)}',$ ${\rm (iii)}',$ ${\rm (iv)}'$ and the
additional one \\ 
${\rm (ii)}''$ 
\[ [\gamma]^2 =Q_{\gamma}, \; \; \;  
{\rm if} \; \; \; \gamma \in \Delta_+, \]  
where $Q_{\gamma}$'s are independent central parameters 
indexed by $\gamma \in \Delta_+.$ The commutative algebra 
generated by 
Dunkl elements in this case may be considered as a 
``multiparameter'' deformation of the coinvariant algebra 
of the Coxeter system $(W,S).$} 
\end{rem}   
\begin{rem} 
{\rm In a similar fashion one can define a 
quantization 
$qBE^+(W,S)$ of the super-version $BE^+(W,S)$ of the bracket 
algebra
$BE(W,S)$ and a family of elements $\tilde{\theta}_s \in qBE^+(W,S).$
See Remark 2.2 for 
the definition of the algebra $BE^+(W,S).$ It is a challenging 
problem to describe the subalgebra in $qBE^+(W,S)$ generated by 
the pairwise anticommutative elements $\tilde{\theta}_s,$ $s\in S.$} 
\end{rem} 
\section{Extended Bruhat graph and quantum Bruhat representation} 

Starting from this section, we assume that the Coxeter system $(W,S)$ is
a crystallographic one. 
Let us denote by $\rho$ the half-sum of all positive roots, i.e.
\[ \rho=\frac{1}{2}\sum_{\gamma\in \Delta_{+}}\gamma. \]
If $\gamma^{\vee}=\sum_{s\in S} n_s \alpha_s^{\vee},$ then 
\[ (\rho,\gamma^{\vee})=\sum_{s\in S} n_s. \] 
\begin{lem} 
Let $\gamma$ be a positive root, then 
\[ 2(\rho,\gamma^{\vee})-1 \geq l(s_{\gamma}). \] 
\end{lem} 
{\it Proof.} For $\gamma \in \Delta_{+},$ define $r$ as the 
minimal number of simple reflections 
$s_0,s_1,\ldots,s_{r-1}$ such that $s_{\gamma}=
s_{r-1}\cdots s_1s_0s_1 \cdots s_{r-1}.$ Then, we can conclude that 
$l(s_{\gamma})=2r-1.$ Hence, by induction on $r$, we have 
$2(\rho,\gamma^{\vee})-1
\geq 2r-1 =l(s_{\gamma}).$ \rule{3mm}{3mm} 
\subsection{Extended Bruhat graph}
\begin{dfn}
The extended Bruhat graph 
$\Gamma(W,S)$ is a graph whose vertices are elements of $W$ with 
arrows $v\rightarrow w$ in the Bruhat ordering and 
additional arrows 
$v\stackrel{\gamma}{\rightarrow}_e w$ which mean that 
$w=vs_{\gamma} \; (\gamma \in \Delta_{+})$ 
and $l(w)=l(v)-2(\rho,\gamma^{\vee})+1.$ 
\end{dfn} 
\begin{lem} 
Let $(W,S)$ be a crystallographic 
Coxeter system and $(W',S')$ its parabolic subsystem. Then the 
extended Bruhat graph $\Gamma(W',S')$ is a subgraph of 
$\Gamma(W,S)$ by the map induced by the inclusion 
$W'\rightarrow W.$ Moreover, if there exists an arrow 
$v\rightarrow_e w$ with $v,w\in W'$ in $\Gamma(W,S),$ 
then the arrow $v\rightarrow_ew$ belongs to $\Gamma(W',S').$ 
\end{lem} 
This follows immediately from Lemma 1.2 and Lemma 7.1.
\begin{rem}
{\rm Definition 7.1 and Lemma 7.2 were discovered originally by
D.Peterson \cite{Pe}. }
\end{rem} 
\subsection{Quantum Bruhat representation}
Let us define an operator $\tilde{\bf s}_{\gamma}$ 
$(\gamma \in \Delta_{+})$ acting on the group ring \\ 
${\bf Q}[q_s | s \in S]
\langle W \rangle ,$ 
by the rule 
\[ \tilde{\bf s}_{\gamma}.w= \left\{ 
\begin{array}{ll} 
ws_{\gamma}, & {\rm if} \; \; l(w)=l(ws_{\gamma})-1, \\ 
q_{\gamma^{\vee}} ws_{\gamma}, & 
{\rm if} \; \; l(w)=l(ws_{\gamma})+2(\rho,\gamma^{\vee})-1, \\ 
0, & {\rm otherwise.}  
\end{array} 
\right.  \] 
\begin{thm} 
A map $[\gamma]\mapsto \tilde{\bf s}_{\gamma}$ defines a 
representation of the quantized bracket algebra $qBE(W,S).$ 
\end{thm} 
{\it Proof}. The compatibility with the relations ${\rm (ii)'}$ 
in Definition 6.1 
is clear. We check the compatibility with the relations 
${\rm (iii)'}$ and ${\rm (iv)'}.$ Let $\Delta'$ be 
as in Definition 2.1 (iii). We are considering only 
crystallographic root systems, so we may assume that 
$\Delta'$ is of type $I_2(m)$ with $m=3,4,6.$ 
Take an arbitrary element $w\in W.$ 
If $l(ws_{\beta}s_{\alpha})=l(w)+2,$ then relation
\[ [\gamma_0][\gamma_{m-1}]w= \sum_i 
[\gamma_i][\gamma_{i+1}]w \] 
follows from the same argument as in the proof of Theorem 3.2. 

Let $\alpha=\gamma_0,$ $\beta=\gamma_{m-1}$ and 
$A(\alpha,\beta)=\{ (\gamma_i,\gamma_{i+1}) | i=0,\ldots,m-2 \} .$
We consider the case $l(ws_{\beta}s_{\alpha})\leq l(w).$ 
Note that $(\rho,\gamma_i^{\vee})\geq (\rho ,\alpha^{\vee})$ 
and 
$(\rho,\gamma_{i+1}^{\vee})\geq (\rho, \beta^{\vee})$ for 
$i=0,\ldots,m-2.$ 
For $(\gamma_1,\gamma_2),(\delta_1,\delta_2) \in 
A(\alpha,\beta),$ 
$(\rho,\gamma_j^{\vee})=(\rho,\delta_j^{\vee})$ holds 
if and only if 
$(\gamma_1,\gamma_2)=(\delta_1,\delta_2).$
Hence, if there exists a path $\Gamma$ of type 
$w\stackrel{\gamma}{\rightarrow}_e
*\stackrel{\delta}{\rightarrow}_e ws_{\beta}s_{\alpha},$ 
then we have 
\[ l(w) = l(ws_{\beta}s_{\alpha})+2(\rho,\gamma^{\vee})
+2(\rho,\delta^{\vee})-2 \geq l(ws_{\beta}s_{\alpha})+
l(s_{\alpha})+l(s_{\beta}) \; \; \; (*). \] 
This means that $l(w)=l(ws_{\beta}s_{\alpha})+l(s_{\alpha})+
l(s_{\beta}),$ and 
$(\gamma,\delta)=(\beta,\alpha).$ In this case, we can see 
that there exists unique pair $(\gamma_1,\gamma_2) \in 
A(\alpha,\beta)$ such that 
if $[\gamma_1][\gamma_2]w \not= 0$ and (*) holds. 
Similarly, if there exists a path $\Gamma$ of type 
$w\stackrel{\beta}{\rightarrow}_e
*\stackrel{\alpha}{\rightarrow} ws_{\beta}s_{\alpha}$ 
or $w\stackrel{\beta}{\rightarrow}
*\stackrel{\alpha}{\rightarrow}_e ws_{\beta}s_{\alpha},$ 
we can find unique pair $(\gamma_1,\gamma_2)\in 
A(\alpha,\beta)$ such that $[\gamma_1][\gamma_2]w \not= 0.$ 
\rule{3mm}{3mm} 
\begin{rem} 
{\rm It follows from our proof that in the extended 
Bruhat graph corresponding to a crystallographic Coxeter group, 
there exist exactly two paths connecting two vertices $v_1,$ 
$v_2$ such that $l(v_1)-l(v_2)\equiv 0$ (mod 2). This property 
does not hold in general for noncrystallographic 
Coxeter systems.} 
\end{rem} 

Now we assume that Coxeter system $(W,S)$ comes from a connected 
simply-connected semi-simple Lie group $G.$ We denote by $B$ 
the Borel subgroup of $G.$ The small quantum cohomology ring 
of the flag variety $G/B$ is isomorphic to the quotient 
ring of the polynomial ring 
${\bf S}(V^{\ast})\otimes {\bf R}[q_s]$ by the ideal 
$\tilde{I}_W$ generated by quantum $W$-invariant polynomials, 
which are explicitly given by Kim \cite{Kim}. 
\begin{thm} 
For Coxeter groups of classical type and 
of type $G_2,$ 
the subalgebra in $qBE(W,S)$ generated by Dunkl elements, 
${\bf R}[q_s][ \tilde{\theta}_s 
\mid s\in S ],$ is canonically isomorphic to the quantum 
cohomology ring $QH^{\ast}(G/B).$ 
\end{thm}  
A proof of Theorem 7.2 is based on  direct computations, see 
Section 9. 
Note that for Lie algebras of type $A,$ Theorem 7.2 was stated 
for the first time in \cite{FK}, and has been proved later 
in \cite{Po}. 
\begin{conj} 
{\rm Theorem 7.2 holds for any crystallographic 
finite Coxeter system.} 
\end{conj} 
\begin{prob} 
{\rm For any finite Coxeter system $(W,S),$ describe ``a quantum 
coinvariant algebra'' of the group $W,$ i.e. to describe
the subalgebra in $qBE(W,S)$ generated by the Dunkl elements 
$\tilde{\theta}_s,$ $s\in S.$} 
\end{prob} 
\section{Quantum Chevalley formula} 
For any polynomial 
$f\in {\bf S}(V^{\ast})\otimes {\bf R}[q_s],$ 
one can define an element $[f]$ of 
$qBE(W,S),$ using the substitution 
$\omega_s \mapsto \tilde{\eta_s}.$ We regard this element 
$[f]$ as an operator acting on the group ring 
${\bf R}[q_s]\langle W \rangle .$ 
\begin{prop}
Let $w\in W,$  there exists a 
unique polynomial 
$\tilde{P}_w\in {\bf S}(V^{\ast})\otimes {\bf R}[q_s]$ 
characterized by the conditions: 
\[ [ \tilde{P}_w](1)=w,\] 
\[ \tilde{P}_w= X_w +\sum_{l(v)<l(w)}c_v X_v \; \; \; \; 
(c_v \in {\bf R}[q_s]), \]
where $X_w$ are the polynomials defined in Section 5, Definition 5.1 .
\end{prop} 
{\it Proof.} 
If $l(w)<2,$ then $[X_w](1)=w.$ In general, 
we have 
\[ [X_w](1)=w+\sum_{l(v)<l(w)} c_v v 
\; \; \; (c_v \in {\bf R}[q_s], \; v \in W). 
\; \; \; \; \rule{3mm}{3mm} \] 
\begin{rem}
{\rm The polynomial $\tilde{P}_w$ defined in 
Proposition 8.1 coincides with the quantum 
Bernstein-Gelfand-Gelfand polynomial 
introduced in \cite{Mar}.} 
\end{rem}

It follows from Theorems 7.1 and 7.2 that for classical Coxeter groups
and $G_2$ one has 
\smallskip \\ 
{\bf Quantum Chevalley formula} (\cite{Pe}, \cite{FW}) \\ 
{\it For $s\in S$ and $w\in W,$ we have 
\[ \tilde{P}_s \tilde{P}_w = 
\sum_{w \stackrel{\gamma}{\rightarrow} w'}
\langle \omega_s,\gamma^{\vee} \rangle \tilde{P}_{w'} 
+ \sum_{w \stackrel{\gamma}{\rightarrow}_e w'}
q_{\gamma^{\vee}} \langle \omega_s,\gamma^{\vee}\rangle 
\tilde{P}_{w'} \; \; {\rm mod} \; \tilde{I}_W, \]
where the sums are taken with respect to the positive roots 
$\gamma.$} \bigskip 

\begin{rem}
{\rm  In Proposition 8.1, we have introduced the 
polynomial $\tilde{P}_w$ 
satisfying the condition $[ \tilde{P}_w](1)=w.$ 
One can consider the action of $[ \tilde{P}_w]$ on 
any element $u\in W$ via the quantum Bruhat representation, 
and obtain an expression 
\[ [ \tilde{P}_w](u)=\sum_{v\in W} c_{wu}^v(q) 
\cdot v, \] 
where $c_{wu}^v(q)\in {\bf R}[q_s]$ are polynomials whose 
coefficients are the so-called $3$-point Gromov-Witten invariants of 
genus zero for the target space $G/B.$ }
\end{rem} 
\begin{conj}
{\rm Let $(W,S)$ be a crystallographic Coxeter system, then there 
exists a monomial basis $\{ b_{\mu} \} _{\mu}$ in the algebra 
$qBE(W,S)$ such that for any $w\in W$ the polynomial 
$[ \tilde{P}_w]$ can be 
written as a linear combination of $b_{\mu}$'s with 
{\it nonnegative} coefficients which do not depend on $q_s$'s. 
}
\end{conj}
\section{Examples} 
Explicit description of relations and the Dunkl elements for 
quantized $A_n$-bracket algebra is given in \cite{FK}. In this Section
we study in more detail the cases of $B_n$-, $D_n$- and $G_2$-bracket
algebras.  \smallskip 

We fix an orthonormal basis $e_1,\ldots ,e_n$ of 
$n$-dimensional 
Euclidean space. 
\subsection{Quantized $B_n$-bracket algebra} 
The root system of type $B_n,$ $n\geq 2,$ consists of 
the elements 
$\pm e_i \pm e_j$ and $\pm e_i$ ($1\leq i,j \leq n$), and 
we fix a set of simple roots 
\[ S(B_n)= \{ \alpha_1=e_1-e_2,\ldots ,\alpha_{n-1}=e_{n-1}-e_n, 
\alpha_n=e_n \} . \] 
The quantized $B_n$-bracket algebra $qBE(B_n)=qBE(W(B_n),S(B_n))$ 
is generated by the symbols
$[i,j]=[e_i-e_j],$ $\overline{[i,j]}=[e_i+e_j]$ 
and $[i]=[e_i]$ 
over ${\bf R}[q_1,\ldots ,q_n]$ subject to the relations: 
\smallskip \\ 
(0) $[i,j]=-[j,i],$ $\overline{[i,j]}=\overline{[j,i]},$ \\ 
(1) $[i,i+1]^2=q_i,$ $[n]^2=q_n,$ \\ 
\hspace*{5.2mm} $[i,j]^2=0$, if $| i-j | \not= 1$;
$[i]^2=0$, if $i<n$; $\overline{[i,j]}^2=0,$ if $i\not= j,$ \smallskip \\ 
(2) $[i,j][k,l]=[k,l][i,j],$ 
$\overline{[i,j]}[k,l]=[k,l]\overline{[i,j]},$ 
$\overline{[i,j]}\overline{[k,l]}=\overline{[k,l]}
\overline{[i,j]},$ \\ 
\hspace*{5.2mm} if $\{ i,j\} \cap \{ k,l\} =\o $, 
\smallskip \\ 
(3) $[i][j]=[j][i],$ $[i,j]\overline{[i,j]}=
\overline{[i,j]}[i,j],$ 
$[i,j][k]=[k][i,j],$ $\overline{[i,j]}[k]=[k]\overline{[i,j]},$ 
if $k\not= i,j$, \smallskip \\ 
(4) $[i,j][j,k]+[j,k][k,i]+[k,i][i,j]=0,$ \\ 
\hspace*{5.2mm} $\overline{[i,k]}[i,j]+[j,i]\overline{[j,k]}
+\overline{[k,j]}\overline{[i,k]}=0,$ \\ 
\hspace*{5.2mm} $[i,j][i]+[j][j,i]+[i]\overline{[i,j]}+
\overline{[i,j]}[j]=0,$ \\
\hspace*{5.2mm} if all $i,$ $j$ and $k$ are distinct,
\smallskip \\  
(5) $[i,j][i]\overline{[i,j]}[i]+\overline{[i,j]}[i][i,j][i]+
[i][i,j][i]\overline{[i,j]}+[i]\overline{[i,j]}[i][i,j]=0,$  if $i<j.$ 
\bigskip \\ 
The Chevalley and Dunkl elements are given by 
$\tilde{\eta}_{s_{\alpha_i}}
=\tilde{\theta}_1+\cdots + \tilde{\theta}_i,$ where 
\[ \tilde{\theta}_i := \tilde{\theta}^{B_n}_i
=\sum_{j\not= i}([i,j]+\overline{[i,j]})+2[i], \; \; \;  1\leq i \leq n. \] 
The Chevalley elements $\tilde{\eta}_{s_{\alpha_i}}$ correspond 
to the Pieri-Chevalley type formula, where as the Dunkl elements 
$\tilde{\theta}_i$ correspond to the Monk type formula in the 
cohomology ring of the flag variety. 
It is easy to see that in the formula for $\tilde{\theta}_i$ above, one can
replace the term $2[i]$ by that $c[i]$ for any constant $c.$ The resulting
operators still commute pairwise.
 
Now we define the quantum $B_n$-invariant polynomials 
following \cite{Kim}. 
Let $E_{i,j}\in M_{2n}({\bf R})$ be a matrix such that 
its $(i,j)$ entry is 1 and other entries are 0. We set 
$t_i=E_{i,i}-E_{i+n,i+n},$ $E_{\alpha_i^{\vee}}=
E_{i+1,i}+E_{i+n,i+n+1},$ 
$E_{-\alpha_i^{\vee}}=E_{i,i+1}-E_{i+n+1,i+n}$ 
($1\leq i \leq n-1$), $E_{\alpha_n^{\vee}}=-2E_{2n,n}$ and 
$E_{-\alpha_n^{\vee}}=2E_{n,2n}.$ Let 
\[ X^B(e,q)=\sum_{i}e_i t_i + \sum_{j}q_j E_{-\alpha_j^{\vee}}+
\sum_{j}E_{\alpha_j^{\vee}}. \] 
The quantum $B_n$-invariant polynomials 
$J^B_{\nu}(e,q)=J^B_{\nu}(e_1,\ldots ,e_n;q_1,\ldots ,q_n)$ 
($1\leq \nu \leq n$) 
are coefficients of the characteristic polynomial of 
$X^B(e,q)$, namely, 
\[ \det (tI +X^B(e,q))=t^{2n}+\sum_{\nu =1}^{n}
J^B_{\nu}(e,q)t^{2(n-\nu )} . \] 
The quantum cohomology ring of $B_n$-flag variety is 
isomorphic to the ring 
\[ {\bf C}[e_1,\ldots, e_n,q_1,\ldots ,q_n]/(J^B_1,\ldots,J^B_n).
\] 
\begin{prop}
In the quantized bracket algebra 
$qBE(B_n)$ we have the following identities 
\[ J^B_{\nu}(\tilde{\theta}_1,\ldots,\tilde{\theta}_n;q)=0, 
\; \; \; 1\leq \nu \leq n. \] 
\end{prop} 
Proof of Proposition 9.1 is based on Lemma 9.1 below. \medskip \\ 
Before to state it, let us introduce a bit of notation. \\ 
{\bf Notation} Let $\{ i,j \}$ denote either generator $[i,j]$ or 
$\overline{[i,j]},$ and define $\overline{ \overline{ [ i ,j ] }}=
[i,j].$ We also define elements $A(a_1,\ldots ,a_k), 
\overline{A}(a_1,\ldots,a_k) \in BE(B_n)$ for distinct integers 
$2\leq a_1,\ldots, a_k \leq n$ as follows 
\[ A(a_1,\ldots,a_k)=\sum_{j=1}^k(-1)^{j-1} \left( 
\prod_{m=j}^k [1, a_m] \right) \cdot [1] \cdot \left( 
\prod_{m=1}^j \overline{ [1, a_m] } \right), \]
\[ \overline{A}(a_1,\ldots,a_k)=\sum_{j=1}^k(-1)^{k-j} \left( 
\prod_{m=j}^k \overline{ [1, a_m] } \right) \cdot [1] \cdot 
\left( \prod_{m=1}^j [1, a_m] \right). \]
\begin{lem} 
We have the following cyclic relations 
in the algebra $BE(B_n)$ for distinct integers 
$2\leq a_1,\ldots, a_k \leq n:$
\smallskip \\ 
$(1)$ $ \{ 1, a_1 \} \{ 1, a_2 \} \cdots 
\{ 1, a_k \} \{ 1, a_1 \} 
+ (${\rm cyclic permutations on indices}$)=0;$  \smallskip \\ 
$(2)$ 
$\{ 1, a_1 \} \{ 1, a_2 \} \cdots \{ 1, a_k \} 
\overline{ \{ 1, a_1 \} } +  \{ 1, a_2 \} \{ 1, a_3 \} 
\cdots \{ 1, a_k \} 
\overline{ \{ 1, a_1 \} }
\overline{ \{ 1, a_2 \} }$ \smallskip \\ 
\hspace*{5mm}
$+ \cdots + \{ 1, a_k \} \overline{ \{ 1, a_1 \} } 
\cdots \overline{ \{ 1, a_{k-1} \} }
\overline{ \{ 1, a_k \} } =  \overline{ \{ 1, a_1 \} } 
\{ 1, a_2 \} 
\cdots \{ 1, a_k \} \{ 1, a_1 \}$ \smallskip \\ 
\hspace*{5mm}
$+ \overline{ \{ 1, a_2 \}} \{ 1, a_3 \} \cdots \{ 1, a_k \} 
\overline{ \{ 1, a_1 \} } \{ 1, a_2 \} + \cdots + 
\overline{ \{ 1, a_k \}} \overline{ \{ 1, a_1 \} } \cdots  
\overline{ \{ 1, a_{k-1} \} } \{ 1, a_k \} ;$ \smallskip \\ 
$(3)$ $[1]( A(a_1,\ldots,a_k)+\overline{A}(a_1,\ldots,a_k) ) + 
( A(a_1,\ldots,a_k)+\overline{A}(a_1,\ldots,a_k) ) [1]=0;$ 
\smallskip \\
$(4)$ All the relations which are obtained from $(1)$, $(2)$ and $(3)$
by the action of the Weyl group. 
\end{lem} 
\begin{ex} 
{\rm For $k=3,$ one can write down the relations in 
Lemma 9.1 
as follows: \smallskip \\
(1) 
$\{ 1,a_1 \} \{ 1, a_2 \} \{ 1, a_3 \} \{ 1, a_1 \} + 
\{ 1,a_2 \} \{ 1, a_3 \} \{ 1, a_1 \} \{ 1, a_2 \} + 
\{ 1,a_3 \} \{ 1, a_1 \} \{ 1, a_2 \} \{ 1, a_3 \}= 0;$ 
\smallskip \\
(2) 
$\{ 1,a_1 \} \{ 1, a_2 \} \{ 1, a_3 \} 
\overline{ \{ 1, a_1 \} } 
+ \{ 1,a_2 \} \{ 1, a_3 \} \overline{ \{ 1, a_1 \} } 
\overline{ \{ 1, a_2 \} } + \{ 1,a_3 \} \overline{ \{ 1, a_1 \} }
\overline{ \{ 1, a_2 \} } \overline{ \{ 1, a_3 \} }$ 
\smallskip \\
\hspace*{5mm} 
$= \overline{ \{ 1,a_1 \} } \{ 1, a_2 \} \{ 1, a_3 \} 
{ \{ 1, a_1 \} } 
+ \overline{ \{ 1,a_2 \} } \{ 1, a_3 \} \overline{ \{ 1, a_1 \} } 
\{ 1, a_2 \} + \overline{ \{ 1,a_3 \} } \overline{ \{ 1, a_1 \} }
\overline{ \{ 1, a_2 \}} \{ 1, a_3 \} ;$ \smallskip \\
(3) $[1][1,a_1][1,a_2][1,a_3][1]\overline{[1,a_1]}-
[1][1,a_2][1,a_3][1]\overline{[1,a_1]}\overline{[1,a_2]}
+[1][1,a_3][1]\overline{[1,a_1]}
\overline{[1,a_2]}\overline{[1,a_3]}$
\smallskip \\ 
\hspace*{5mm} 
$+[1]\overline{[1,a_1]}\overline{[1,a_2]}
\overline{[1,a_3]}[1][1,a_1] 
-[1]\overline{[1,a_2]}\overline{[1,a_3]}[1][1,a_1][1,a_2]+
[1]\overline{[1,a_3]}[1][1,a_1][1,a_2][1,a_3]$ 
\smallskip \\ 
\hspace*{5mm} 
$+[1,a_1][1,a_2][1,a_3][1]\overline{[1,a_1]}[1]-
[1,a_2][1,a_3][1]\overline{[1,a_1]}\overline{[1,a_2]}[1]
+[1,a_3][1]\overline{[1,a_1]}\overline{[1,a_2]}
\overline{[1,a_3]}[1]$
\smallskip \\ 
\hspace*{5mm}
$+\overline{[1,a_1]}\overline{[1,a_2]}
\overline{[1,a_3]}[1][1,a_1][1] 
-\overline{[1,a_2]}\overline{[1,a_3]}[1][1,a_1][1,a_2][1]+
\overline{[1,a_3]}[1][1,a_1][1,a_2][1,a_3][1]$ 
\smallskip \\ 
\hspace*{5mm}
$=0.$ }
\end{ex} 
In the final part of this Subsection we consider an application of
Lemma 2.1 to the case of $B_n$-bracket algebra. \\
Let $x_i = [i,n],$ $y_i=\overline{[i,n]}$ for $1\leq i \leq n-1,$ 
and $z_n=[n]$ be elements of $BE(B_n).$ Denote by 
${\cal A}_n^B={\cal A}(\Delta(B_n) \setminus \Delta(B_{n-1}))$ 
the subalgebra of $BE(B_n)$ generated by 
$x_1,\ldots , x_n, y_1, \ldots ,y_n$ and $z_n.$
\begin{prop} Action of the twisted derivation $D_{[\gamma]}$ on the 
algebra ${\cal A}^B_n$ is determined 
by the following formulas: 
\[ 
D_{[i,j]}(x_i)=-x_ix_j,\;\;  D_{[i,j]}(x_j)=x_jx_i, \; \; 
D_{[i,j]}(y_i)=-y_iy_j, \; \;  D_{[i,j]}(y_j)=y_jy_i, \]
\[ D_{\overline{[i,j]}}(x_i)=x_iy_j, \; \; 
D_{\overline{[i,j]}}(x_j)=x_jy_i, \; \;  
D_{\overline{[i,j]}}(y_i)=y_ix_j, \; \;  
D_{\overline{[i,j]}}(y_j)=y_jx_i, \]
\[ D_{[i]}(x_i)=x_iz_n-z_ny_i, \; \; D_{[i]}(y_i)=z_nx_i-y_iz_n, \]
\[ D_{[i,j]}(z_n)=D_{\overline{[i,j]}}(z_n)=D_{[i]}(z_n)=0, \] 
for $1\leq i \not = j \leq n-1,$ and 
\[ D_{ \{ k,l \} }(x_i)=D_{ \{ k,l \} }(y_i) = 
D_{[k]}(x_i)=D_{[k]}(y_i)= 0, \; \; 
\textrm{if i, k, l are all distinct,} \]
and the twisted Leibniz rule, see definition of the former in Section 2.2.
\end{prop} 
Therefore, the subalgebra ${\cal A}^B_n$ 
is invariant under the twisted derivation $D_{[\gamma]}$ 
for any root $\gamma\in \Delta(B_{n-1}).$ 
By applying Lemma 2.1 successively, we obtain the following 
decomposition of the algebra $BE(B_n)$ for $n\geq 2,$ 
\[ BE(B_n) \cong {\cal A}^B_2 \otimes 
\cdots \otimes {\cal A}^B_n. \] 
Note that the relations in Lemma 9.1 can be obtained by applying the twisted 
derivations successively to the defining relations 
of the bracket algebra. We don't know whether or not all the relations 
in the bracket algebra can be obtained in such a way. 
\begin{ex}
{\rm By applying $D_{[a_2,a_3]}D_{[a_1,a_2]}$ to the 
4-term relation 
\[ [1][1,a_1][1]\overline{[1,a_1]}+
[1,a_1][1]\overline{[1,a_1]}[1]+[1]\overline{[1,a_1]}[1][1,a_1]
+\overline{[1,a_1]}[1][1,a_1][1]=0, \] 
we obtain the relation (3) in Example 9.1.}
\end{ex}
We conclude this subsection by a construction of one more 
representation of the algebra $BE(B_{n-1}).$ 
Denote by ${\cal F}^B_n$ the quotient of the free associative 
algebra over ${\bf R}$ generated by $X_1,\ldots ,X_{n-1}, 
Y_1,\ldots ,Y_{n-1}$ and $Z_n$ modulo the two-side ideal generated by 
$X_i^2,$ $Y_i^2,$ $Z_n^2$ and 
$Z_n X_i Z_n Y_i + X_i Z_n Y_i Z_n + Z_n Y_i Z_n X_i + 
Y_i Z_n X_i Z_n$ for $1\leq i \leq n-1.$ The Weyl group $W(B_{n-1})$ 
acts on the algebra ${\cal F}^B_n$ by the rule 
\[ s_{ij}(X_i)= X_j, \; s_{ij}(Y_i)=Y_j, \; 
s_{\overline{ij}}(X_i)=-Y_j, \; s_{\overline{ij}}(Y_i) = -X_j , \; 
s_{i}(X_i) = -Y_i, \; s_i(Y_i)=-X_i, \]
\[ s_{ij}(X_k)=s_{\overline{ij}}(X_k)=s_{i}(X_k)= X_k, \; 
s_{ij}(Y_k)=s_{\overline{ij}}(Y_k)=s_{i}(Y_k)=Y_k, \]
\[ s_{ij}(Z_n)=s_{\overline{ij}}(Z_n)=s_i(Z_n)=Z_n \] 
for distinct $i,j,k \in \{ 1,\ldots ,n-1 \} .$  
Now define operators $\nabla_{[i,j]},$ 
$\nabla_{\overline{[i,j]}}$ and $\nabla_{[i]},$ 
$1\leq i \not= j \leq n-1,$ which act on the algebra 
${\cal F}^B_n$ 
by the same formulas as for the operators $D_{[i,j]},$ 
$D_{\overline{[i,j]}}$ and $D_{[i]}$ 
from Proposition 9.2 after replacing $x_i,$ $y_i$ and $z_n$ 
by $X_i,$ $Y_i$ and $Z_n$ respectively. 
Then the operators $\nabla_{[i,j]}$ and 
$\nabla_{\overline{[i,j]}}$ and $\nabla_{[i]},$ 
$1\leq i \not= j \leq n-1,$ give rise to a representation of 
the algebra $BE(B_{n-1})$ in the algebra ${\cal F}^B_n,$ and 
natural epimorphism $\pi ^B_n:{\cal F}^B_n \rightarrow 
{\cal A}^B_n$ is compatible with the action of the algebra 
$BE(B_{n-1}).$ 
\begin{prob}
{\rm Describe the kernel of the epimorphism $\pi ^B_n.$}
\end{prob}
\subsection{Pieri formula for $B_n$-bracket algebra}
The main goal of this subsection is to describe a 
$B_n$-analog of Pieri's formula in some cases, namely, 
we give an explicit 
formula for the value of elementary symmetric polynomials 
of arbitrary degree and 
complete symmetric polynomials of degree two 
in the bracket algebra $BE(B_n)$ after the substitution 
of variables by the $B_n$-Dunkl elements. Let us observe that if 
we specialize all the generators $[i] \in BE(B_n)$ to zero, we 
obtain a $D_n$-analog of Pieri's formula. To state our result, 
it is convenient to introduce a bit of notation. 
Let $S= \{ i_1 <i_2 < \cdots <i_s :=r \}$ be a set of positive 
integers. Define inductively a family of elements 
$\{ K_l(S) \}_{l \geq 1}$ in the algebra $BE(B_r)$ by the following 
rules \smallskip \\ 
{\abovedisplayshortskip -4mm 
i)
\[ K_1(S)= \sum_{i\in S} [i] + \sum_{i \leq j ,\;  i,j \in S} 
\overline{[i,j]} \; ; \; K_l(S)=0, \; \; \; \textrm{if} \; \; \; 
s<l; \] 
ii) 
\[ K_l(S) = K_l(S \setminus \{ r \} ) + \sum_{a\in S} 
([a,r]+\overline{[a,r]})K_{l-1}(S\setminus \{ a \} ) 
+ K_{l-1}(S\setminus \{ r \} )\theta_{r,S}, \] }
where $\theta_{r,S}= \sum_{a\in S}(-[a,r]+\overline{[a,r]})+2[r].$ 
\begin{thm} 
{\rm (1a)} Let $m\leq n,$ then 
\[ e_k(\theta^{B_n}_1,\ldots , \theta_m^{B_n}) = 
\tilde{\sum_{(*)}} \prod_{a=1}^k \{ i_a,j_a \} + 
2 \sum_{l=1}^k \tilde{\sum_{(*)}} \prod_{a=1}^{k-l} \{ i_a ,j_a \} 
K_l(\{ 1,\ldots ,m \} \setminus \{ i_1, \ldots ,i_{k-l} \}), 
 \] 
where the symbol $\tilde{\sum }$ means that in the corresponding sums 
we have to take only distinct monomials among the products 
$\prod_{a=1}^k \{ i_a,j_a \}$ and $\prod_{a=1}^{k-l} \{ i_a,j_a \} ;$ 
the condition $(*)$ means that $1\leq i_a \leq m <j_a 
\leq n$ and all indices $i_a$ are distinct. 
\smallskip \\
{\rm (1b)} The elements $K_l(S)$ can be expressed in the algebra 
$BE(B_r)$ as a linear combination of monomials in $[i]$'s and 
$\overline{[ i,j ]}$'s with nonnegative integer coefficients. 
\smallskip \\
{\rm (1c)} If the number of elements in the set $S$ is equal 
to $l,$ then $K_l(S)=0$ after the specialization $[a]=0$ for all $a \in S.$ 
\smallskip \\ 
{\rm (1d)} $K_2(S)=(K_1(S))^2;$ 
\[ K_3(S) = K_3( S\setminus \{ r \} )+ \sum_{a\in S} 
\overline{[a,r]} K_2(S\setminus \{ a \} ) +K_2(S\setminus \{ r \} ) 
\left( \sum_{a\in S} \overline{[a,r]}+2 [r] \right) \] 
\[ + \sum_{a\in S} \left( \overline{[a,r]}[a] + [r] \overline{[a,r]} + 
\sum_{b\in S} \overline{[b,r]} \overline{[a,b]} \right) 
K_1(S\setminus \{ a \} ) + K_1(S\setminus \{ r \} ) 
\sum_{a\in S} \left( \overline{[a,r]}[a] + [r] \overline{[a,r]} + 
\sum_{b\in S} \overline{[b,r]} \overline{[a,b]} \right) . \] 
For example, the multiplicity of the monomial $\overline{[12]}
\overline{[34]}\overline{[56]}$ in $e_3(\theta_1^{B_6}, 
\theta_2^{B_6},\ldots , \theta_6^{B_6})$ is equal to 4. \\
{\rm (2)}  Let $m\leq n,$ then 
\[ h_2(\theta^{B_n}_1,\ldots , \theta_m^{B_n}) = 
\tilde{\sum_{(**)}}  \prod_{a=1}^2 \{ i_a,j_a \} + 
2 \sum_{1\leq i \leq m <j \leq n} \{ i ,j \} 
K_1(\{ 1,\ldots ,m \}\setminus \{ i \} ) +2 K_2(\{ 1,\ldots , m \} ) \]
\[ +2 \tilde{\sum_{1\leq i_a \leq m < j \leq n}} \left( 
\overline{[i_1,j]}[i_2,j]+[i_1,j] \overline{[i_2,j]} \right) 
+ 2 \sum_{1 \leq i \leq m < j \leq n} 
\left( \{ i,j \} [i] + [i] \{ i,j \} \right) , \]
where the condition $(**)$ means that $1\leq i_a \leq m <j_a 
\leq n$ and all $j_a$'s are distinct. 
\smallskip \\
{\rm (3)} 
$h_k(\theta_1^{B_n},\ldots, \theta_m^{B_n})=0,$ if 
$k+m>2n.$ 
\end{thm} 
Finally, let us remark that for classical Coxeter 
groups $W=W(A_n),$ $W(B_n),$ and $W(D_n),$ the condition 
$|R(u)|=1,$ $u\in W,$ is equivalent to the condition that 
modulo the ideal $I_W,$ the Schubert class $X_u$ is equal to 
either $e_k({\bf X}_m)$ or $h_k({\bf X}_m)$ for some 
$k$ and $m\leq n,$ up to multiplication by some power of 2. 
In the case of symmetric groups, the permutations $w$
such that $|R(w)|=1$ are exactly the permutations of the following forms: \\
$w= h(a,b):= ( 1,2,...,a+b,a,b+1,b+2,... ),$  or \\
$w= e(a,b):= (1,2,...,a-b,a-b+2,...,a+1,a-b+1,a+2,a+3,... ),$ see e.g. 
\cite{Sot}.
\subsection{Quantized $B_2$-algebra and quantum cohomology}
Here we give an explicit calculation of 
quantum cohomology ring and certain polynomial representatives 
for Schubert classes for $B_2$-flag variety.  
The bracket algebra $qBE(B_2)$ is generated by 
the symbols $[12],$ $\overline{[12]},$ $[1]$ and $[2]$ subject 
to the following relations: \\ 
(i) $[12]^2=q_1,$ $\overline{[12]}^2=0,$ 
$[1]^2=0,$ $[2]^2=q_2,$ \smallskip \\ 
(ii) $[12]\overline{[12]}=\overline{[12]}[12],$ 
$[1][2]=[2][1],$ \smallskip \\ 
(iii) $ [12][1]-[2][12]+[1]\overline{[12]}+\overline{[12]}[2]=0,$ 
$[1][12]-[12][2]+\overline{[12]}[1]+[2]\overline{[12]}=0,$ 
\smallskip \\ 
(iv) $[12][1]\overline{[12]}[1]+\overline{[12]}[1][12][1]+
[1][12][1]\overline{[12]}+[1]\overline{[12]}[1][12]=0.$ 
\medskip \\ 
The Chevalley and Dunkl elements are 
$\tilde{\eta}_{s_{\alpha_1}}=
\tilde{\theta}_1$ and 
$\tilde{\eta}_{s_{\alpha_2}}=\tilde{\theta}_1+\tilde{\theta}_2,$ 
where 
\[ \tilde{\theta}_1=[12]+\overline{[12]}+2[1], \; \; 
\tilde{\theta}_2=-[12]+\overline{[12]}+2[2]. \] 
Quantum cohomology ring of the $B_2$-flag 
variety is isomorphic to the algebra 
${\bf C}[q_1,q_2][e_1,e_2]/I_{B_2},$ 
where 
\[ I_{B_2}=\left(
e_1^2+e_2^2-2q_1-4q_2,e_1^2e_2^2+2q_1e_1e_2-4q_2e_1^2+q_1^2
\right) . \] 
The subalgebra generated by $\tilde{\theta}_1, 
\tilde{\theta}_2$ in 
$qBE(B_2)\otimes {\bf C}$ is isomorphic to the quantum 
cohomology ring. 
Let us consider the quantum Bruhat representation of 
$qBE(B_2)$ and regard the Dunkl elements $\tilde{\theta}_i$ 
as operators acting on the group ring 
${\bf R}[q_1,\ldots,q_n]\langle W(B_2) \rangle.$ Denote by 
$s_{12}$ and $s_2$ the simple reflections with respect to 
the simple roots $e_1-e_2,$ and $e_2$ respectively. 
Then, 
\[ \frac{\tilde{\theta}_1^2-q_1}{2}(id.)=s_{2}s_{12} \] 
\[ \frac{\tilde{\theta}_1\tilde{\theta}_2+q_1}{2}(id.)=
s_{12}s_{2} \] 
\[ \frac{\tilde{\theta}_1^3-2q_1\tilde{\theta}_1-q_1
\tilde{\theta}_2}{2}(id.)
=s_{12}s_{2}s_{12} \] 
\[ \frac{\tilde{\theta}_1^2\tilde{\theta}_2-
\tilde{\theta}_1^3+3q_1\tilde{\theta}_1+q_1\tilde{\theta}_2}{4}
(id.) 
=s_{2}s_{12}s_{2} \] 
\[ \frac{\tilde{\theta}_1^3\tilde{\theta}_2+
q_1\tilde{\theta}_1^2-q_1\tilde{\theta}_1\tilde{\theta}_2 
-q_1^2-4q_1q_2}{4}(id.)=s_{12}s_{2}s_{12}s_{2}. \] 
\begin{rem} 
{\rm Both algebras $BE(B_2)$ and $BE^+(B_2)$ are infinite 
dimensional, but if we add the new relation 
 \[  [1][12][1][12]=[12][1][12][1]  \]
in the algebra $BE(B_2),$ and that
\[  [1][12][1][12]+[12][1][12][1]=0  \]
in the algebra $BE^+(B_2),$ the resulting algebras appear to be finite
 dimensional and have the same Hilbert polynomial
\[ (1+t)^4 (1+t^2)^2 . \]
One can check that the pointed Hopf algebra over the Coxeter group 
$D_4$ constructed in \cite{MS}, is isomorphic to the quotient of the 
algebra $BE^+(B_2)$ by the relation of degree $4$ defined above.} 
\end{rem} 
\subsection{Quantized $D_n$-bracket algebra} 
In $D_n$ case, $n\geq 2,$ fix a set of simple roots as 
\[ S(D_n)=\left\{ \alpha_1=e_1-e_2,\ldots ,
\alpha_{n-1}=e_{n-1}-e_n,\alpha_n=e_{n-1}+e_n \right\} . \] 
The quantized $D_n$-bracket algebra $qBE(D_n)=qBE(W(D_n),S(D_n))$ 
is generated by the symbols
$[i,j]=[e_i-e_j]$ and $\overline{[i,j]}=[e_i+e_j]$ 
over ${\bf R}[q_1,\ldots ,q_n]$ subject to the following 
relations: 
\smallskip \\ 
(0) $[i,j]=-[j,i],$ $\overline{[i,j]}=\overline{[j,i]},$ \\ 
(1) $[i,i+1]^2=q_i,$ $\overline{[n-1,n]}^2=q_n,$ \\ 
\hspace*{5.2mm} $[i,j]^2=0,$ if $|i-j| \not= 1$, 
$\overline{[i,j]}^2=0,$ if $(i,j)\not=(n-1,n),(n,n-1)$, 
\smallskip \\ 
(2) $[i,j][k,l]=[k,l][i,j],$ 
$\overline{[i,j]}[k,l]=[k,l]\overline{[i,j]},$ 
$\overline{[i,j]}\overline{[k,l]}=
\overline{[k,l]}\overline{[i,j]},$ \\ 
\hspace*{5.2mm} if  $\{ i,j\} \cap \{ k,l\} =\o $, 
\smallskip \\ 
(3) $[i,j]\overline{[i,j]}=\overline{[i,j]}[i,j],$ 
\smallskip \\ 
(4) $[i,j][j,k]+[j,k][k,i]+[k,i][i,j]=0,$ 
$\overline{[i,k]}[i,j]+[j,i]\overline{[j,k]}
+\overline{[k,j]}\overline{[i,k]}=0,$ \\
\hspace*{5.2mm} if all $i,$ $j$ and $k$ are distinct. 
\begin{rem}
{\rm Our construction of the quantized bracket algebra $qBE(D_n)$
is compatible with the isomorphisms between the 
Coxeter systems $D_2 \cong A_1 \times A_1$ and $D_3 \cong 
A_3.$ It is easy to see that $qBE(D_2)\cong qBE(A_1) 
\times qBE(A_1)$ and $qBE(D_3) \cong qBE(A_3).$ }
\end{rem}
We set 
$t_i=E_{i,i}-E_{i+n,i+n},$ $E_{\alpha_i^{\vee}}=-E_{i+1,i}+
E_{i+n,i+n+1},$ 
$E_{-\alpha_i^{\vee}}=E_{i,i+1}-E_{i+n+1,i+n}$ 
($1\leq i \leq n-1$), $E_{\alpha_n^{\vee}}=-E_{2n-1,n}+
E_{2n,n-1}$ and 
$E_{-\alpha_n^{\vee}}=E_{n,2n-1}-E_{n-1,2n}.$ Let 
\[ X(e,q)=\sum_{i}e_i t_i + \sum_{j}q_j E_{-\alpha_j^{\vee}}+
\sum_{j}E_{\alpha_j^{\vee}}. \] 
We define the polynomials 
$J^D_{\nu}(e,q)$ by the equation 
\[ \det (tI +X^D(e,q))=t^{2n}+\sum_{\nu =1}^{n}
J^D_{\nu}(e,q)t^{2(n-\nu )} . \] 
Then the quantum cohomology ring of $D_n$-flag variety 
is isomorphic to the ring 
\[ {\bf C}[e_1,\ldots, e_n,q_1,\ldots ,q_n]/
(J^D_1,\ldots,J^D_{n-1},\overline{J^D_n}), \] 
where $\overline{J^D_n}$ is a polynomial such that 
$(\overline{J^D_n})^2=J^D_n.$ 

The Chevalley and Dunkl elements are given 
by $\tilde{\eta}_{s_{\alpha_i}}
=\tilde{\theta}_1+\cdots + \tilde{\theta}_i,$ where 
\[ \tilde{\theta}_i =\sum_{j\not= i}([i,j]+\overline{[i,j]}),
\; \; \; 1\leq i \leq n. \] 
\begin{prop} 
In the quantized bracket algebra 
$qBE(D_n)$ we have the following identities 
\[ J^D_{\nu}(\tilde{\theta}_1,\ldots,\tilde{\theta}_n;q)=0, 
\; \; \; 1\leq \nu \leq n. \] 
\end{prop} 
Proof of Proposition 9.2 follows from the following lemma. 
\begin{lem} 
Relations $(1),$ $(2)$ in Lemma 9.1 and all the relations 
obtained from them by the action of the Weyl group, 
hold also in the algebra $BE(D_n).$ 
\end{lem} 
\begin{rem}
{\rm The non-quantized bracket algebra $BE(D_n)$ is 
a quotient ring of $BE(B_n)$ obtained by putting $[i]=0,$ 
and the Dunkl elements of $BE(D_n)$ are images of those 
of $BE(B_n).$ Hence the Dunkl elements of $BE(D_n)$ 
satisfy the equations coming from the $B_n$ case. 
However, $qBE(D_n)$ is not a quotient of $qBE(B_n).$} 
\end{rem} 
\begin{ex}
{\rm  Quantum $D_n$-invariants for $n=4,$ \medskip \\ 
$ J_1^D=  - e_1^2- e_2^2- e_3^2- e_4^2+2 q_1 +2 q_2 +2 q_3 +2 q_4, $ 
\bigskip \\ 
$ J_2^D=  q_3^2 +2 q_1 e_1 e_2+ q_4^2-2 q_4 e_3 e_4 +2 q_2 q_4 +2 q_3 e_3 
e_4 + q_2^2-2 q_3 q_4 +2 q_1 q_2$ \medskip \\ 
$+4 q_1 q_3 +2 q_2 q_3 + q_1^2 +4 q_1 q_4 
+2 q_2 e_2 e_3 -2 q_1 e_3^2-2 q_2 e_4^2-2 q_4 e_1^2+ e_1^2 
e_3^2-2 q_2 e_1^2 $ \medskip \\ 
$ + e_3^2 e_4^2 + e_1^2 e_2^2-2 q_1 e_4^2 -2 q_3  e_1^2+ 
e_1^2 e_4^2+ e_2^2 e_4 ^2 + e_2^2 e_3^2-2 q_4 e_2^2-2 q_3 e_2^2, $ 
\bigskip \\ 
$J_3^D=  -2 q_2 q_3 e_1^2- e_1^2 e_3^2 e_4^2-2 q_2 q_4 e_1^2+2 
q_4 e_1^2 e_2^2-2 q_1 q_2 e_4 ^2 - e_2^2 e_3^2 e_4^2+2 q_3 e_1 ^2 e_2^2 $ 
\medskip \\ 
$+ 2 q_1 e_3^2 e_4^2- e_1^2 e_2^2 e_4^2+4 q_1 q_3 e_3 e_4 -4 q_1 q_3 q_4 
+2 q_1 q_3^2 +2 q_2 e_1^2 e_4^2- e_1^2 e_2^2 e_3^2+2 q_3 q_4 e_1^2 $ 
\medskip \\ 
$+2 q_1 q_2 q_3 +2 q_1 q_2 q_4+2 q_3 q_4 e_2^2- q_1^2 e_4^2- q_1^2 e_3^2
+2 q_4 q_1^2- q_4^2 e_1^2+2 q_4^2 q_1 +2 q_3 q_1^2 $ \medskip \\ 
$- q_3^2 e_1^2 - q_3^2 e_2^2- q_2^2 e_1^2- q_2 ^2 e_4^2-
2 q_1 q_2 e_1 e_3-2 q_2 q_3 e_2 e_4 -2 q_1 e_1 e_2 e_3^2 +
4 q_1 q_4 e_1 e_2 $ \medskip \\ 
$+2 q_4 e_1^2 e_3 e_4-2 q_2 e_1^2 e_2 e_3 
+4 q_1 q_3 e_1 e_2 -4 q_1 q_4 e_3 e_4- q_4^2 e_2^2-2 q_3 e_2^2 e_3 e_4 $ 
\medskip \\ 
$ -2 q_1 e_1 e_2 e_4^2 +2 q_4 e_2^2 e_3 e_4-
2 q_3 e_1^2 e_3 e_4+2 q_2 q_4 e_2 e_4 -2 q_2 e_2 e_3 e_4^2 , $ \bigskip \\ 
$\overline{J_4^D}= e_1 e_2 e_3 e_4 + q_1 e_3 e_4 + q_2 e_1 e_4 
+ q_3 e_1 e_2 - q_4 e_1 e_2 + q_1 q_3 -q_1 q_4.$ }
\end{ex}
\begin{rem}
{\rm We don't know whether or not the algebra $BE(D_4)$ is finite 
dimensional. However, the commutative quotient of the algebra 
$BE(D_4)$ is finite dimensional and has the following Hilbert 
polynomial 
\[ 1+12t+50t^2+84t^3+48t^4=(1+2t)(1+4t)(1+6t+6t^2)=
(1+t)(1+3t)^2 (1+5t)+3t^4. \] 
Let us remark that the polynomial $(1+t)(1+3t)^2 (1+5t)$ coincides 
with the Hilbert polynomial of the cohomology ring of the pure 
braid group of type $D_4.$ 

It was a big surprise for us to find 
that the Hilbert polynomial of the commutative quotient of 
the algebra $BE(D_5)$ is equal to 
\[ 1+20t+150t^2+520t^3+824t^4+480t^5=(1+2t)(1+4t)(1+6t)(1+8t+10t^2). \]
However, the obvious generalization of the above formulas for the Hilbert 
polynomial of the commutative quotient of the algebra $BE(D_n)$ 
is false. 
}
\end{rem}
Similar to the case of $B_n$-bracket algebra, 
the subalgebra ${\cal A}^D_n={\cal A}(\Delta(D_n) 
\setminus \Delta(D_{n-1}))$ 
generated by $x_i=[i,n]$ and $y_i=\overline{[i,n]},$ 
$i=1,\ldots,n-1,$ in the algebra $BE(D_n)$ 
is invariant under the twisted derivation $D_{[\gamma]}$ 
for any root $\gamma \in \Delta (D_{n-1}).$ By applying Lemma 2.1
successively, we obtain the following decomposition 
of the algebra $BE(D_n)$ for $n\geq 2,$ 
\[ BE(D_n) \cong {\cal A}^D_2\otimes \cdots \otimes 
{\cal A}^D_n. \] 
\subsection{Quantized $G_2$-bracket algebra}
Fix a set of positive roots of type $G_2$ as 
\[ \{ a, b=3a+f, c=2a+f,d=3a+2f,e=a+f,f \} . \] 
Then quantized $G_2$-bracket algebra is generated by the symbols 
$a,b,c,d,e,f$ with the relations: \\ 
(1) $a^2=q_1,$ $f^2=q_2,$ $b^2=c^2=d^2=e^2=0,$ \\
(2) $ea=ce+ac,$ $ae=ec+ca,$ $fb=df+bd,$ $bf=fd+db,$ \\ 
\hspace*{5mm} $eb=be,$ $cf=fc,$ $ad=da,$ \\ 
\hspace*{5mm} $af=ba+cb+dc+ed+fe,$ $fa=ab+bc+cd+de+ef,$ \\ 
(3) $bcdefd+dbcdef+fedcbd+dfedcb=0,$ \\
\hspace*{5mm} $fabcdb+bfabcd+dcbafb+bdcbaf=0,$ \\
\hspace*{5mm} $defabf+fdefab+bafedf+fbafed=0.$ \\
The Chevalley elements are defined by 
\[ \tilde{\eta}_{s_a}=a+3b+2c+3d+e, \] 
\[ \tilde{\eta}_{s_f}=b+c+2d+e+f. \] 
Let $\tilde{\theta}_1=\tilde{\eta}_{s_a}-
\tilde{\eta}_{s_f}$ and $\tilde{\theta}_2=\tilde{\eta}_{s_f}$ 
be the corresponding Dunkl elements, 
then we have the relations $g_2(\tilde{\theta}_1,\tilde{\theta}_2)=
g_6(\tilde{\theta}_1,\tilde{\theta}_2)=0$ 
in the algebra $qBE(G_2),$ where 
\[g_2(\xi_1,\xi_2):=\xi_1^2+\xi_2^2-\xi_1 \xi_2-q_1-3q_2, \]  
\[g_6(\xi_1,\xi_2):= \xi_1^3\xi_2^3-3q_2\xi_1^2\xi_2^2+q_1\xi_1\xi_2^3+
q_1\xi_2^4+q_1q_2\xi_1^2+3q_2(q_1+q_2)\xi_1\xi_2+
2q_1q_2\xi_2^2+q_1^2q_2-6q_1q_2^2-q_2^3. \] \\
The small quantum cohomology ring of $G_2$-flag variety is isomorphic to 
the ring 
\[ {\bf R}[q_1,q_2][\xi_1,\xi_2]/(g_2,g_6). \]  
One can check the latter representation for the small quantum 
cohomology ring of $G_2$-flag variety is equivalent to that given by
B.Kim \cite{Kim}.

\subsection{Dunkl elements and fundamental invariant polynomials for $I_2(m)$}
Let $a_i= \mu_i e_1 + \lambda_i e_2,$ where $\mu_i = \cos(i \pi /m)$ 
and $\lambda_i = \sin(i \pi /m)$ for $i=0,1,\ldots ,m-1.$ 
Then $\Delta_+ = \{ a_1 ,\ldots , a_{m-1} \}$ forms the set of 
positive roots of type $I_2(m).$ The set of simple roots is 
$S=\{ a_0,a_{m-1} \}$ and 
\[ a_i = \lambda_1^{-1}\lambda_{i+1}a_0+ \lambda_1^{-1}
\lambda_i a_{m-1}. \] 
The Chevalley elements 
are given by 
\[ \eta_{s_{a_1}}= \sum_{i=0}^{m-1}\lambda_1^{-1}\lambda_{i+1}[a_i], \] 
\[ \eta_{s_{a_{m-1}}} = \sum_{i=0}^{m-1}\lambda_1^{-1}\lambda_i[a_i]. \] 
Let 
$\theta_1 = \eta_{s_{a_0}} + \lambda_1 \eta_{s_{a_{m-1}}}$ 
and 
$\theta_2 = (\lambda_1^{-1} \mu_1 +1)\eta_{s_{a_0}} + 
(\lambda_1^{-1} + \mu_1)\eta_{s_{a_{m-1}}}$ be the Dunkl 
elements of type $I_2(m).$ 
The fundamental invariant polynomials are 
\[ f_2(\xi_1,\xi_2)= \xi_1 ^2 + \xi_2 ^2, \] 
and 
\[ f_m(\xi_1,\xi_2) = \sum_{i=0}^{[m/2]} (-1)^i {m \choose 2i}
\xi_1^{2i}\xi_2^{m-2i} , \]
where $\xi_1$ and $\xi_2$ are variables corresponding to the 
orthonormal basis $e_1$ and $e_2.$ 
\begin{prop}
In the algebra $BE(I_2(m))$ one has 
\[ f_2(\theta_1 ,\theta_2)= 0, \; \; f_m(\theta_1,\theta_2)=0. \]
\end{prop}
We can check that the algebra generated by the Dunkl elements 
$\theta_1$ and $\theta_2$ in the algebra $BE(I_2(m))$ is 
isomorphic to the quotient of the polynomial 
ring ${\bf R}[\xi_1 ,\xi_2 ]/(f_2,f_m).$ 
\begin{rem}
{\rm In this subsection, all roots are normalized to satisfy 
the condition $(a_i,a_i)=1.$ The root systems of type $I_2(4)$ and 
$I_2(6)$ can be identified with the crystallographic systems of 
type $B_2$ and $G_2,$ but 
the choice of the normalization is different. Hence, the 
Dunkl elements for $I_2(4)$ and $I_2(6)$ in this subsection have 
a different expression from the ones defined in Subsections 9.3 and 9.5. }
\end{rem}

Research Institute for Mathematical Sciences \\ 
Kyoto University \\
Sakyo-ku, Kyoto 606-8502, Japan \\
e-mail: kirillov@kurims.kyoto-u.ac.jp 
\bigskip \\ 
Department of Mathematics \\ 
Kyoto University \\ 
Sakyo-ku, Kyoto 606-8502, Japan \\ 
e-mail: maeno@math.kyoto-u.ac.jp 

\begin{thebibliography}{95}
\bibitem{AS} N. Andruskiewitsch and H.-J. Schneider, {\it 
Pointed Hopf algebras,} New directions in Hopf algebras, 
Math. Sci. Res. Inst. Publ. {\bf 43} Cambridge Univ. Press, 
Cambridge, 2002, 1-68.
\bibitem{BS} N. Bergeron and F. Sottile, {\it A Pieri-type 
formula for isotropic flag manifolds,} Trans. Amer. Math. Soc. 
{\bf 354} (2002), 2659-2705. 
\bibitem{BGG} I. N. Bernstein, I. M. Gelfand and S. I. Gelfand, 
{\it Schubert cells and the cohomology of the spaces $G/P,$} 
Russian Math. Surveys {\bf 28} (1973), 1-26. 
\bibitem{De} M. Demazure, {\it Invariants sym\'etriques entiers 
des groupes de Weyl et torsion,}  \\
Inv. Math. {\bf 21} (1973), 
287-301. 
\bibitem{DFI} P. Di Francesco and C. Itzykson, {\it Quantum 
intersection rings,} in {\it The Moduli Space of Curves,} 
(R. Dijkgraaf, C. Faber and G. van der Geer, eds.) Progress in 
Math. {\bf 129} Birkh\"auser, 1995, 81-148. 
\bibitem{Du} C. F. Dunkl, {\it Harmonic polynomials and peak sets 
of reflection groups,} Geom. Dedicata 
{\bf 32} (1989), 157-171. 
\bibitem{FGP} S. Fomin, S. Gelfand and A. Postnikov, {\it 
Quantum Schubert polynomials,} Journ. Amer. Math. Soc. {\bf 10} 
(1997), 565-596.
\bibitem{FK} S. Fomin and A. N. Kirillov, {\it Quadratic algebras, 
Dunkl elements and Schubert calculus,} in {\it Advances in 
Geometry,} (J.-L. Brylinski, R. Brylinski, V. Nistor, B. Tsygan 
and P. Xu, eds.) Progress in Math. {\bf 172} Birkh\"auser, 
1995, 147-182. 
\bibitem{FP} S. Fomin and C. Procesi, {\it Fibered quadratic 
Hopf algebras related to Schubert calculus,}  \\
J. Algebra {\bf 230} 
(2000), 174-183. 
\bibitem{FW} W. Fulton and C. Woodward, {\it On the quantum 
product of Schubert classes,} math.AG/0112183. 
\bibitem{GK} A. Givental and B. Kim, {\it Quantum cohomology 
of flag manifolds and Toda lattices,}  \\
Commun. Math. Phys. {\bf 168} (1995), 609-641. 
\bibitem{Hi}H. L. Hiller, {\it Schubert calculus of a Coxeter 
group,} Enseign. Math. {\bf 27} (1981), 57-84. 
\bibitem{Hu}J. E. Humphreys, {\it Reflection Groups and 
Coxeter Groups,} Cambridge Studies in Adv. Math. {\bf 29} 
Cambridge University Press, 1990. 
\bibitem{Kim} B. Kim, {\it Quantum cohomology of flag 
manifolds $G/B$ and quantum Toda lattices,}  \\
Annals of Math. {\bf 149} (1999), 129-148. 
\bibitem{Kir} A. N. Kirillov, {\it On some quadratic 
algebras,} L. D. Faddeev's Seminar on Mathematical Physics, 
Amer. Math. Soc. Transl. Ser. 2 {\bf 201}, AMS, 
Providence, RI, 2000, 91-113. 
\bibitem{KM} A. N. Kirillov and T. Maeno, {\it Quantum 
double Schubert polynomials, quantum Schubert polynomials 
and Vafa-Intriligator formula,} Discrete Math. {\bf 217} 
(2000), 191-223. 
\bibitem{LS} A. Lascoux and M.-P. Sch\"utzenberger, {\it 
Polyn\^omes de Schubert,} C. R. Acad. Sci. Paris {\bf 294} 
(1982), 447-450. 
\bibitem{Maj} S. Majid, {\it Noncommutative differentials and 
Yang-Mills on permutation groups $S_N,$} math/0105253. 
\bibitem{Mar} A.-L. Mare, {\it The combinatorial quantum 
cohomology ring of $G/B,$} math.CO/0301257. 
\bibitem{MS} A. Milinski and H.-J. Schneider, {\it Pointed 
indecomposable Hopf algebras over Coxeter groups,} Contemp. 
Math. {\bf 267} (2000), 215-236. 
\bibitem{Pe} D. Peterson, {\it Quantum cohomology 
of flag varieties, Lectures given at MIT,} 1997. 
\bibitem{Po} A. Postnikov, {\it On a quantum version of Pieri's 
formula,} in {\it Advances in Geometry,}  \\ 
(J.-L. Brylinski, R. Brylinski, V. Nistor, B. Tsygan 
and P. Xu, eds.) Progress in Math. {\bf 172} Birkh\"auser, 
1995, 371-383. 
\bibitem{PR} P. Pragacz and J. Ratajski, {\it Formulas 
for Lagrangian and orthogonal degeneracy loci; $\tilde Q$-polynomial 
approach,} Compositio Math. {\bf 107} (1997), 11-87. 
\bibitem{Sot} F. Sottile, {\it Pieri's rule for flag manifolds and
Schubert polynomials,}  \\
Annales Fourier {\bf 46} (1996), 89-110. 
\bibitem{Wo} S. L. Woronowicz, {\it Differential calculus 
on compact matrix pseudogroups (quantum groups),} Commun. Math. 
Phys. {\bf 122} (1989), 125-170. 
\end{thebibliography}
\end{document}